\documentclass[reqno]{amsart}
\newtheorem{thm}{Theorem}[section]
\newtheorem{df}[thm]{Definition}
\newtheorem{prop}[thm]{Proposition}

\newtheorem{lem}[thm]{Lemma}
\newtheorem{rem}[thm]{Remark}

\newtheorem{cor}[thm]{Corollary}
\newtheorem{exm}[thm]{Example}

\usepackage{amsmath,amssymb,amscd,epsfig}
\usepackage{color, graphics}
\usepackage{epsfig,texdraw}
\usepackage{latexsym}
\usepackage{mathrsfs}
\usepackage{eucal}
\usepackage{hyperref}
\numberwithin{equation}{section}
\newcommand{\Ld}{\Lambda}
\newcommand{\ld}{\lambda}

\newcommand{\pf}{\noindent{\bfseries Proof. }}
\renewcommand{\baselinestretch}{1.2}
\begin{document}

\begin{abstract}
The purpose of this this paper is to generalize the functors arising from the theory of Witt vectors duto to Cartier.
Given a polynomial $g(q)\in \mathbb Z[q]$,
we construct a functor ${\overline {W}}^{g(q)}$ from the category of $\mathbb Z[q]$-algebras to that of commutative rings.
When $q$ is specialized into an integer $m$, it produces a functor from the category of commutative rings with unity to that of commutative rings.
In a similar way, we also construct several functors related to ${\overline { W}}^{g(q)}$.
Functorial and structural properties such as induction, restriction, classification and unitalness
will be investigated intensively.
\end{abstract}
\title[Generalizing Witt vector construction]
{Generalizing Witt vector construction}

\author{YOUNG-TAK OH}
\address{Department of Mathematics \\ Sogang University\\Seoul 121-742, Korea}
\email{ytoh@sogang.ac.kr}
\maketitle
\baselineskip=12pt

\section{Introduction}
\renewcommand{\thefootnote}{}
\footnotetext[1]{%
\renewcommand{\baselinestretch}{1.2}\selectfont
\noindent
This research was supported by Basic Science Research Program
through the National Research Foundation of Korea(NRF) funded
by the Ministry of Education, Science and Technology
(2012R1A1A2001635).
\hfill \break 2010 MSC :{13F35 (primary), 19A22, 20E18 (secondary)}
\hfill \break Keywords: Witt vector, Witt ring, necklace ring, $\ld$-ring,
Adams operation, symmetric function.}

The big Witt ring functor $W$ of Cartier \cite{Ca} was introduced as a generalization of the classical $p$-typical Witt
ring functor $W_p$ of Teichm\"{u}ller and Witt \cite{W}.
Since its birth, it had remained as a mysterious object for a long time due to its bizarre structure.
For instance, even their value at
the ring of integers had been veiled until the beginning of nineteen-eighties.

The necklace ring was contrived by Metropolis and Rota \cite{MR}
to understand Witt vectors and their operations in a combinatorial way.
Given a commutative ring $A$ with unity,
they endowed the set of infinite sequences with entries in $A$ with a commutative ring structure
such that the addition is defined componentwise and
the multiplication is defined as
$$(a_n)_{n\in \mathbb N}\cdot (b_n)_{n\in \mathbb N}=(c_n)_{n\in \mathbb N}, \text{ where }c_n=\sum_{[i,j]=n}(i,j)a_jb_j.$$
Here $[i,j]$ denotes the least common multiple of $i$ and $j$,
$(i,j)$ the greatest common divisor of $i$ and $j$, and $\mathbb N$
the set of positive integers, respectively.
This construction induces a functor $Nr$ and the ring $Nr(A)$ is called {\it the necklace ring over $A$}.
It turns out that $W(A)$ is functorially isomorphic to $Nr(A)$ over binomial rings
(see \cite {O2006}).

The main objective of this paper is to generalize the construction of $W$ and $Nr$.
Let {\textsf RINGS} be the category of commutative rings (not necessarily unital) and ring homomorphisms,
{\textsf Rings} the category of commutative rings with unity and unity-preserving ring homomorphisms, and
{\textsf {$\mathbb Z[q]$-Algebras}} the category of $\mathbb Z[q]$-algebras with unity and
unity-preserving $\mathbb Z[q]$-algebra homomorphisms. Here $q$ is an indeterminate.
We construct a functor
$${\overline { W}}^{g(q)}: \textsf{$\mathbb Z[q]$-Algebras} \to \textsf{RINGS}$$
attached to each polynomial $g(q)\in \mathbb Z[q]$.
When $g(q)=q$ and $q$ is specialized into an integer $m$, it produces a functor
$$W^{m}: \textsf{Rings} \to \textsf{RINGS}.$$
This process can be illustrated as follows:
\begin{equation*}
\begin{CD}
{\overline { W}}^{q}: \textsf{$\mathbb Z[q]$-Algebras} \to \textsf{RINGS}\\
@VV{q \mapsto m}V  \\
W^{m}: \textsf{Rings} \to \textsf{RINGS}
\end{CD}
\end{equation*}
Viewed as a functor from \textsf{Rings} into itself, $W^{0}$ coincides with $W$.

For categories $\mathcal C_i \,\,(i=1,2)$,
let $F^{\mathcal C_1}_{\mathcal C_2}:\mathcal C_1 \to \mathcal C_2$ denote
the forgetful functor assigning to each object $R$ in $\mathcal C_1$ the object $R$ in $\mathcal C_2$.
In the context of Witt vectors, it is well known that
$W$ is isomorphic to the $\ld$-ring functor $\Lambda$ due to Grothendieck.
To be more precise, we have
$$ W\cong F^{\textsf{$\lambda$-Rings}}_{\textsf {Rings}} \circ \Lambda.$$
(see Section \ref{Witt ring functor and necklace ring functor}).
This isomorphism naturally led us to construct functor ${\overline \Lambda}^{g(q)}$ and ${\Lambda}^{m}$ such that
\begin{equation*}
{\rm (1)}\,\,{\overline \Lambda}^{g(q)}\cong {\overline {W}}^{g(q)} \text{ and } {\rm (2)}\,\,\Lambda^m\cong W^{m}.
\end{equation*}

It is quite interesting to notice that the relationship between $W$ and $Nr$ can be extended to our setting.
To do this, we construct a functor
$${\overline {\textsf B}}^{g(q)}:{\textsf{$\mathbb Z[q]\otimes {\Psi}$-Rings}} \to \textsf{RINGS},$$
where {\textsf{$\mathbb Z[q]\otimes {\Psi}$-Rings}} is the category of $\mathbb Z[q]$-algebras equipped with
the $\Psi$-ring structure determined by $\Psi^n(q)=q^n \,\,(n\ge 1).$
The morphisms are $\mathbb Z[q]$-algebra homomorphisms compatible with $\Psi$-operations.
We also construct $$\textsf B^{m}:{\textsf{${\Psi}$-Rings}} \to \textsf{RINGS}$$
for each integer $m$.
It should be noted that $\textsf B^{m}$ cannot be obtained directly from ${\overline {\textsf B}}^{q}$ by specializing $q$ into $m$.
It turns out that these functors have close connections to ${\overline {W}}^{g(q)}$ and $W^m$.
Indeed we show that
$${\overline {W}}^{g(q)} \circ F^{{\textsf{$\mathbb Z[q]\otimes {\lambda}$-Rings}}}_{\textsf{Rings}}
\cong {{{\overline {\textsf B}}}}^{g(q)}|_{{\textsf{$\mathbb Z[q]\otimes {\lambda}$-Rings}}}.$$
and
$$W^m \circ F^{\textsf{$\ld$-Rings}}_{\textsf{Rings}} \cong {\textsf B}^m|_{\textsf {$\ld$-Rings} }.$$
Here {\textsf{$\mathbb Z[q]\otimes {\lambda}$-Rings}} is the subcategory of
{\textsf{$\mathbb Z[q]\otimes {\Psi}$-Rings}}
whose objects are $\mathbb Z[q]$-algebras equipped with the $\ld$-ring structure determined by $\psi^n=\Psi^n \,\,(n\ge 1)$
and whose morphisms are $\mathbb Z[q]$-algebra homomorphisms preserving $\ld$-operations.

This paper is organized as follows:
In Section \ref{Preliminaries}, we introduce prerequisites necessary to develop our arguments.
In Section \ref{$q$-deformation of Witt rings and its specializations}, we construct ${\overline {W}}^{g(q)}$
for each $g(q)\in \mathbb Z[q]$ and $W^{m}$ for each integer $m\in \mathbb Z$.
We also classify ${{{\overline {W}}}}^{g(q)}$ and
$F^{\textsf{RINGS}}_{\textsf {Abel}}\circ{{{\overline {W}}}}^{g(q)}$
up to strict-isomorphism as $g(q)$ ranges over $\mathbb Z[q]$.
The classification of ${{W}}^{m}$ and
$F^{\textsf{RINGS}}_{\textsf {Abel}}\circ W^m$
up to strict-isomorphism
will also be provided as $m$ ranges over the set of integers.
In Section \ref {q-Deformation of Grothendieck's Lambda functor},
we provide functors ${\overline \Lambda}^{g(q)}$ and $\Lambda^{m}$ which are equivalent to ${\overline {W}}^{g(q)}$
and $W^{m}$, respectively.
Some new bases of $\mathbb Q[q]\otimes \Lambda_{\mathbb Z}$ and $\mathbb Z[q, \frac {1}{1-q}]\otimes \Lambda_{\mathbb Z}$
will be also given, where $\Lambda_{\mathbb Z}$ is the ring of symmetric functions in infinitely many variables with integer coefficients.
In Section \ref{q-Deformations of necklace rings and their specializations}, we construct
${\overline {\textsf B}}^{g(q)}$ and $\textsf B^{m}$.
In particular, in Sec \ref{Structure constants},
structure constants of ${\overline {\textsf B}}^{g(q)}$ and ${\textsf B}^m$ will be intensively investigated.
These are particularly interesting due to their connection to {\it M\"{o}bius-like functions}.
In Section \ref{Functorial and structural properties}, we study functorial and structural properties of the functors introduced in the previous sections.
To be more precise, we deal with induction, restriction, and the conditions for the unitalness.
In Section \ref{Relationship between Witt rings and necklace rings}, we demonstrate the connection between
${\overline {W}}^{g(q)}$ and ${{\overline {\textsf B}}}^{g(q)}$ and that of $W^{m}$ and $\textsf B^m.$
As an application, we give a ring theoretical interpretation of the product identities recently introduced by Kim and Lee \cite{KL}
(see Remark \ref{Henry Kim}).

\section{Preliminaries}
\label{Preliminaries}
Here we collect definitions, terms, and notation which will be required to develop our arguments.
\subsection{Pre-$\ld$-rings, $\ld$-rings, and $\Psi$-rings}
For the clarity we will use the character $R$ for pre-$\ld$-rings, $\ld$-rings, and $\Psi$-rings.
On the other hand, the character $A$ will be used to denote a commutative ring.

\subsubsection{}
\label{subsection on pre lambda rings}
A {\it pre-$\lambda$-ring} $R$ is a commutative ring with unity, ${1}_R$ (or simply $1$), together with
a series of operations
$\lambda^n : R \to R, \,
(n=0,1,2,\ldots),$ such that, for all $x,y \in R$,
\begin{align*}
(1)\,\,\lambda^0(x)={1}_R, \quad
(2)\,\,\lambda^1(x)=x, \quad
(3)\,\,\lambda^n(x+y)=\sum_{r=0}^{n}\lambda^r(x)\lambda^{n-r}(y).
\end{align*}
Let $R$ and $R'$ be pre-$\lambda$-rings.
A ring homomorphism $f:R\to R'$ with $f({1}_R)={1}_{R'}$ is
called a {\it pre-$\lambda$-ring homomorphism} if
$\lambda^n(f(x))=f(\lambda^n(x))$ for all $x\in R$ and all $n\ge 0$.

Given a positive integer $n$, we define the {\it $n$th Adams operation, $\psi^n :R\to R$},
so that $(-1)^{n-1}\psi^{n}(x)$ is the coefficient of $\frac{d}{d{\it t}}\log
\lambda_t(x)$, that is,
\begin{equation}\label {how to define Adams operations}
\frac{d}{d{\it t}}\log
\lambda_t(x)=\sum_{n=1}^{\infty}(-1)^{n-1}\psi^{n}(x)t^{n-1}, \quad \forall x\in R.
\end{equation}
It is easy to show that
$\psi^1$ is the identity function and
$\psi^n$ is additive for all $n\ge 1$.
Comparing the coefficient of either side of Eq.\eqref{how to define Adams operations} yields
the well known \textit{Newton's formula}:
\begin{equation*}
\underset{i=0}{\stackrel{n-1}{\sum}}(-1)^i \psi^{n-i}(x) \ld^i(x)  + (-1)^nn\ld^n(x)= 0,
\quad \forall \,x \in R, \,\forall \,n\ge 1.
\end{equation*}
Utilizing this formula repeatedly,
one can deduce that $\psi^n(x)$ is a polynomial in $\ld^i(x)$ ($1 \le i \le n$)
with integer coefficients.
To be more precise, for each positive integer $n$,
there exists the unique polynomial $\mathcal Q_n\in \mathbb Z[t_1,\ldots, t_n]$
with $\psi^n(x) =\mathcal Q_n(\ld^1(x),\ld^2(x),\ldots,\ld^n(x)).$
Hence, every pre-$\lambda$-ring homomorphism
preserves Adams operations.

\subsubsection{}
\label{Review on Lambda rings}
Naively speaking, a $\ld$-ring $R$ is a pre-$\lambda$-ring
with very special properties for $\ld^n(xy)$ and
$\ld^n(\ld^m(x))$ for all $x,y\in R$ and $m,n \ge 1$.

Let $X=\{x_i:i\ge 1 \}$ and $Y=\{y_i: i\ge 1
\}$ be the infinite sets of commuting indeterminates $x_i$ and
$y_i$ $(i\ge 1)$, respectively. We call these sets {\it alphabets}.
Let us introduce the following operations on alphabets:
\begin{align*}
&X+Y=\{x_i,y_i: i\ge 1\}, \quad X\cdot Y=\{x_iy_j: i,j\ge 1\}, \text{ and }\\
&\Lambda^m(X)=\{x_{i_1}\cdots x_{i_m}: i_1<i_2<\cdots <i_m\} \quad (m \ge 1).
\end{align*}
Let $h_n(X)$ be the $n$th complete symmetric function in $x_i$ $(i\ge 1)$.

\begin{lem}{\rm (\cite{AT})} With the above notation, we have the following.

{\rm (a)} For each positive integer $n$,
$h_n(XY)\in \mathbb Z[h_1(X),\ldots,h_n(X);h_1(Y),\ldots, h_n(Y)].$

{\rm (b)} For all positive integers $m$ and $n$,
$h_n(\Lambda^m(X))\in \mathbb Z[h_i(X): 1\le i \le mn].$
\end{lem}

Let $P_n\in \mathbb Z[s_1,\ldots,s_n;t_1,\ldots,t_n]$ be the unique
polynomial characterized by
$$P_n(h_1(X),\ldots, h_n(X);h_1(Y), \ldots, h_n(Y))=h_n(XY)$$
and $P_{n,m}\in \mathbb Z[s_i: 1\le i \le mn]$ the unique polynomial
characterized by
$$P_{n,m}(h_1(X),\ldots, h_{mn}(X))=h_n(\Lambda^m(X)).$$

\begin{df}\label{def of lambda ring}
A pre-$\ld$-ring $R$ is called a $\ld$-ring if the following conditions are satisfied:
\begin{enumerate}
\item $\ld_t({1}_R) ={1}_R+t$ .
\item For all $x,y  \in  R$ and $n \geq 0$,
$\ld^n(xy)=P_n(\ld^1(x),\ldots,\ld^n(x);\ld^1(y),\ldots,\ld^n(y)).$
\item
For all $x\in R$ and $m,n \geq 0$, $\ld^n(\ld^m(x))=P_{n,m}(\ld^1(x),\ldots,\ld^{nm}(x)).$
\end{enumerate}
\end{df}
For more information on $\ld$-rings, refer to \cite{AT,H1,H2,DK,O2005, O2006,Y}.

\subsubsection{}
A {\it $\Psi$-ring} $R$ is a
commutative ring with unity together with a series of operations $\Psi^n: R\to
R\,\,\, (n \ge 1)$ such that, for all $x,y\in R$ and $m,n\ge 1$,
\begin{enumerate}
\item
$ \Psi^1(x)=x, $
\item
$\Psi^n(x+y)=\Psi^n(x)+\Psi^n(y),$
\item
$\Psi^n({1}_R)={ 1}_R,$
\item
$\Psi^n(xy)=\Psi^n(x)\Psi^n(y),$ and
\item
$\Psi^n(\Psi^m(x))=\Psi^{nm}(x).$
\end{enumerate}
Let $R,R'$ be $\Psi$-rings.
A ring homomorphism $f:R\to R'$ with $f({1}_R)={1}_{R'}$ is called a {\it $\Psi$-ring homomorphism}
if $\Psi^n(f(x))=f(\Psi^n(x))$ for all $x\in R$ and all $n\ge 1$.

Note that $\Psi$-operations can be interpreted as Adams operations of a $\ld$-ring in some cases.
\begin{prop}{\rm (\cite{W})}\label {psi=ld}
Let $R$ be a $\Psi$-ring which has no $\mathbb Z$-torsion
and such that $\Psi^p(a)=a^p \pmod {pR}$ if $p$ is a prime.
Then there is the unique $\lambda$-ring structure on $R$ such that
the $\Psi^n\, (n\ge 1)$ are the associated Adams operations.
\end{prop}

\subsection{Categories and forgetful functors}
\label{Categories and forgetful functors}
In the present paper, we deal with many kinds of ring categories and forgetful functors.
Here we collect all of them for the understanding of the readers.
\begin{enumerate}

\item Let {\textsf Abel} be the category of abelian groups and group homomorphisms.

\item Let {\textsf RINGS} be the category of commutative rings (not necessarily unital) and ring homomorphisms.

\item Let {\textsf Rings} be the category of commutative rings with unity and unity-preserving ring homomorphisms.

\item Let {\textsf {$\mathbb Z[q]$-Algebras}} be the category of $\mathbb Z[q]$-algebras with unity and
unity-preserving $\mathbb Z[q]$-algebra homomorphisms. Here $q$ is an indeterminate.

\item Let \textsf {Pre-$\lambda$-Rings} be the category of pre-$\lambda$-rings and pre-$\lambda$-ring homomorphisms.

\item Let \textsf{$\lambda$-Rings} be the category of $\lambda$-rings and pre-$\lambda$-ring homomorphisms.

\item Let \textsf{$\Psi$-Rings} be the category of $\Psi$-rings and $\Psi$-ring homomorphisms.

\item Let {\textsf{$\mathbb Z[q]\otimes {\Psi}$-Rings}} be the category of $\mathbb Z[q]$-algebras that are equipped with
a $\Psi$-ring structure with $\Psi^n(q)=q^n$ for all $n\ge 1.$
The morphisms are $\mathbb Z[q]$-algebra homomorphisms compatible with $\Psi$-operations.

\item
Let {\textsf{$\mathbb Z[q]\otimes {\lambda}$-Rings}} be the subcategory of
{\textsf{$\mathbb Z[q]\otimes {\Psi}$-Rings}}
such that the objects are $\mathbb Z[q]$-algebras equipped with the $\ld$-ring structure with $\psi^n=\Psi^n$ for all $n\ge 1$
and the morphisms are $\mathbb Z[q]$-algebra homomorphisms preserving $\ld$-operations.
\end{enumerate}

Let $F^{\mathcal C_1}_{\mathcal C_2}:\mathcal C_1 \to \mathcal C_2$ be
the forgetful functor assigning to each object $R$ in $\mathcal C_1$ the object $R$ in $\mathcal C_2$.
For instance, we will use the following forgetful functors.
\begin{enumerate}
\item Let
$F^{\textsf{$\ld$-Rings}}_{\textsf{Rings}}:\textsf{{$\lambda$-Rings}}\to \textsf{Rings}$ be the forgetful functor assigning to each
$\lambda$-ring $R$ its underlying ring $R$.

\item Let
$F^{\textsf{$\Psi$-Rings}}_{\textsf{Rings}}:\textsf{{$\Psi$-Rings}}\to \textsf{Rings}$ be
the forgetful functor assigning to each $\Psi$-ring $R$ its underlying ring $R$.

\item Let
$F^{\textsf{RINGS}}_{\textsf {Abel}}:\textsf{RINGS} \to \textsf{Abel}$ be the forgetful functor assigning to each
commutative ring $A$ its underlying abelian group $A$.

\item Let
$F^{\textsf{$\mathbb Z[q]$-Algebras}}_{\textsf{Rings}}:\textsf{$\mathbb Z[q]$-Algebras} \to \textsf{RINGS}$ be the forgetful functor assigning to
each $\mathbb Z[q]$-algebra $A$ its underlying ring $A$.
\end{enumerate}

Finally, when $\mathcal C$ is a subcategory of $\mathcal D$, the notation $\iota^{\mathcal C}_{\mathcal D}:\mathcal C \to \mathcal D$
denotes the natural injection induced by the assignment $A\to A$.
On the other hand, for a functor $\mathscr F: D \to \mathcal D'$, the notation $\mathscr F|_{\mathcal C}: \mathcal C \to \mathcal D'$
denotes the restriction of $\mathscr F$ onto $\mathcal C$.

\subsection{Witt ring functor and necklace ring functor}
\label{Witt ring functor and necklace ring functor}
Finally, let us review basic definitions and notation on the ring-valued functors
$W, \Lambda$, and $B$ very briefly.
Let ${gh}:\textsf{Rings} \to \textsf{Rings}$ be the unique functor defined by
the following conditions:
\begin{enumerate}
\item
As a set,
${gh}(A)=A^\mathbb N=\{(a_n)_{n\in \mathbb N}: a_n\in A \text { for all }n\ge 1\}$.
\item
Addition and multiplication are defined componentwise.
\item
For every unity-preserving ring homomorphism $f:A\to B$
and every $\alpha \in {gh}(A)$ one has
${gh}(f)(\alpha)=f \circ \alpha$.
\end{enumerate}
The following characterization of the big Witt
ring functor is due to Cartier.

\begin{prop}\label{definition of witt ring functor}{\rm (\cite{Ca})}
There exists the unique functor
$$W:\textsf{Rings} \to \textsf{Rings}$$
satisfying the following conditions$:$
\begin{enumerate}
\item
As a set,
$ W(A)=A^\mathbb N$.
\item
For every unity-preserving ring homomorphism $f:A\to B$
and every $\alpha \in W(A)$ one has $W(f)(\alpha)=f\circ \alpha$.
\item
The map $\Phi_A:  W(A)\to {gh}(A)$ defined by
$$(a_n)_{n\in \mathbb N} \mapsto
\left(\sum_{d|n}da_d^{\frac nd}\right)_{n\in \mathbb N}$$
is a ring homomorphism.
\end{enumerate}
\end{prop}

Next, we introduce Grothendieck's $\Ld$-functor.
More information can be found in \cite{AT,H1,H2,DK,O2005,O2006}.
Given a commutative ring $A$ with unity, let
$\Lambda(A)$ be the set of all formal power series of the form
${1} + \sum_{n=1}^{\infty}a_nt^n, \,a_n \in A \text{ for all }n\ge 1.$
Grothendieck showed that $\Lambda(A)$ can have a
$\ld$-ring structure with respect to the operations defined in the following manner:
\begin{enumerate}
\item
Addition, $+_{\Ld}$, is the usual product of power series.

\item
Multiplication, $\times_{\Ld}$, is given by
\begin{equation*}
\left({1}  + \sum_{n\ge 1}a_nt^n \right)\times_{\Ld} \left(1  + \sum_{n\ge 1}b_nt^n \right)
={1}+ \sum_{n\ge 1}P_n(a_1,\dots,a_n;b_1,\dots,b_n)t^n \ .
\end{equation*}

\item
For each nonnegative integer $m$,
the $m$th $\ld$-operation, $\lambda^m$, is given by
\begin{equation*}
\lambda^m \left(1 + \sum_{n\ge 1}a_nt^n \right)
=1 + \sum_{n=1}P_{n,m}(a_1,\dots,a_{nm})t^n \ .
\end{equation*}
\end{enumerate}
For the definition of $P_n$ and $P_{n,m}$, refer to Section \ref{subsection on pre lambda rings}.
Given a unity-preserving ring homomorphism $f:A\to B$, define
$\Lambda(f):\Lambda(A) \to \Lambda(B)$ by
$$\Lambda(f)\left(1+\sum_{n\ge 1}a_nt^n \right)=1+\sum_{n\ge 1}f(a_n)t^n.$$
One can easily sees that the assignment $(A, f) \mapsto (\Ld(A), \Ld(f))$ induces a functor $\Lambda:\textsf{Rings} \to \textsf{$\ld$-Rings}$.

The functor $\Lambda$ has a characterization similar to Proposition \ref{definition of witt ring functor}.
Let ${Gh}:\textsf{Rings} \to \textsf{Rings}$ be the functor defined by
the following conditions:
\begin{enumerate}
\item
As a set,
${Gh}(A)=A[[t]]=\{\sum_{n\ge 1}a_nt^n: a_n \in A\}$.
\item
Addition and multiplication are defined componentwise.
\item
For every unity-preserving ring homomorphism $f:A\to B$
and every $\alpha \in {Gh}(A)$ one has
$Gh(f)(\alpha)=f \circ \alpha$.
\end{enumerate}
With this preparation, we can state the following characterization.
\begin{prop} {\rm (\cite{H1,H2,DK})}
The functor $\Lambda$ is uniquely determined by the following conditions$:$
\begin{enumerate}
\item
As a set, $\Lambda(A)=\{{1}+\sum_{n\ge 1}a_nt^n : a_n \in A\}$.
\item
For every unity-preserving ring homomorphism $f:A\to B$, one has
\begin{align*}
\Lambda (f)
:\Lambda (A)
\to \Lambda (B),\,{1}+\sum_{n\ge 1}a_nt^n \mapsto {1}+\sum_{n\ge 1}f(a_n)t^n.
\end{align*}
\item
The map $\frac{d}{dt}\log:  \Lambda  (A) \to Gh(A)$ defined by
$${1}+\sum_{n\ge 1}a_nt^n \mapsto
\frac{\sum_{n\ge 1}na_nt^{n-1}}{1+\sum_{n\ge 1}a_nt^n}$$
is a ring homomorphism.
\end{enumerate}
\end{prop}

We have a natural isomorphism $\Theta: W \to F^{\textsf{$\ld$-Rings}}_{\textsf{Rings}}\circ \Lambda$, where
$$\Theta_A: W(A) \to  \Lambda(A) \quad (a_n)_{n\in \mathbb N} \mapsto \prod_{n\ge 1}\frac{1}{1-a_nt^n}.$$
For more information, refer to \cite{H1}.

In the rest of this subsection, we introduce the functor $$B:\textsf{$\Psi$-Rings} \to \textsf{Rings}.$$
As stated in the introduction, Metropolis and Rota \cite{MR} introduced the necklace ring functor $Nr$ such that
${Nr}(\mathbb Z)$ is isomorphic of $W(\mathbb Z)$.
In \cite{O2005}, ${Nr}$ was generalized to
$$Nr_G: \textsf{$\ld$-Rings} \to \textsf{Rings},$$
where $G$ is an arbitrary profinite group.
In particular, if $G=\hat C$, the profinite completion of the multiplicative infinite cyclic group $C$,
and $R$ is a $\ld$-ring, then $Nr_{\hat C}(R)$ turns out to be the ring $R^{\mathbb N}$
such that the addition is defined componentwise and
the multiplication is defined by
$$(a_n)_{n\in \mathbb N}\cdot (b_n)_{n\in \mathbb N}=\left(\sum_{[i,j]=n}(i,j)\psi^{\frac ni}(a_i)\psi^{\frac nj}(b_j)\right)_{n\in \mathbb N}.$$
The subsequent proposition tells us that $Nr_{\hat C}$ can be extended to the category of $\Psi$-rings.

\begin{prop} {\rm (\cite{O:2012})}
There exists the unique functor
$$B:\textsf{$\Psi$-Rings} \to \textsf{Rings}$$
satisfying the following conditions:
\begin{enumerate}
\item
As a set,
$B(R)=R^\mathbb N$.
\item
For every $\Psi$-ring homomorphism $f:R\to S$
and every $\alpha \in B(R)$ one has
$B(f)(\alpha)=f \circ \alpha$.
\item
The map $\varphi_{R}: {B}(R)\to {gh}(R)$ defined by
$$(a_n)_{n\in \mathbb N} \mapsto
\left(\sum_{d|n}d\Psi^{\frac nd}(a_d)\right)_{n\in \mathbb N}$$
is a ring homomorphism.
\end{enumerate}
\end{prop}

\begin{lem}\label{to define Teichmuller map}{\rm (\cite[Section 3]{O2006})}
Let $x$ be an indeterminate and
$E$ the $\ld$-ring freely generated by $x$.
For every positive integer $n$, let
$$M(x,n):=\frac 1n \sum_{d|n}\mu(d)\psi^d(x^{\frac nd}).$$
Then $M(x,n)\in \mathbb Z[\ld^i(x):1\le i\le n]$.
\end{lem}

Given a $\ld$-ring $R$ and an element $r\in R$,
let $\pi: \textsf E \to R$ be the unique $\ld$-ring homomorphism
with $\pi(x)=r$
and let $M(r,n)$ be $\pi(M(x,n))$ for all $n\ge 1$.
Given a $\ld$-ring $R$, consider the map
$$\tau_R: W(R)\to B(R),
\quad (a_n)_{n\in \mathbb N}\mapsto \left(\sum_{d|n}M(a_d,n/d)\right)_{n\in \mathbb N}.
$$
The following theorem is about the relationship between $W$ and $B$.

\begin{thm} {\rm (\cite[1.9.3]{E1} and \cite{O2005,O2006})}\label{Theorem of connection}
The assignment $R\mapsto \tau_R$ induces a natural isomorphism $\tau:  W\circ
F^{\textsf{$\ld$-Rings}}_{\textsf{Rings}}\to B|_{\textsf{$\ld$-Rings}}$.
\end{thm}

\section{$q$-deformation of Witt rings and its specializations}
\label{$q$-deformation of Witt rings and its specializations}

Hereafter $q$ denotes an indeterminate and $g(q)$ any polynomial in $\mathbb Z[q]$.
In this section, we construct two functors
$${\overline { W}}^{g(q)},\,{\overline \Lambda}^{g(q)}:\textsf{$\mathbb Z[q]$-Algebras} \to \textsf{RINGS}$$
which are naturally isomorphic to each other.
When $g(q)=q$, we obtain functors
$$W^m,\,\Lambda^m:\textsf{Rings} \to \textsf{RINGS}$$
by specializing $q$ into an arbitrary integer $m$.

\begin{df}

{\rm (a)} Let $p$ be a prime.
For any nonzero integer $a$, we define $\nu_p(a)$ to be the $p$-valuation of $a$,
the exponent of $p$ in the highest power of $p$ dividing $a$.

{\rm (b)}
Let $q$ be an indeterminate and $m$ an integer.
For each positive integer $n$, we let
$[n]_q:=1+q+\cdots+q^{n-1}$ and
$[n]_m:=1+m+\cdots+m^{n-1}.$
\end{df}

The following lemma, which plays a crucial role in our arguments, can be easily verified by using the binomial expansion formula.
\begin{lem}\label{key modular equivalence 1}
Let $n$ be a positive integer and $p$ any prime divisor of $n$.

{\rm (a)}
For any positive integer $k$ and any polynomial $f(x_1,\ldots,x_k)\in \mathbb Z[x_1,\ldots,x_k]$,
it holds that
$f(x_1^p,\ldots,x_k^p)^{\frac np}\equiv f(x_1,\ldots,x_k)^{n}\pmod {p^{\nu_p(n)}}.$

{\rm (b)} Let $m$ be an integer.
Given a divisor $d$ of $n$, let $p$ be a prime divisor of $\frac nd$. Then
$d(1-m^{\frac nd}) \equiv d(1-m^{\frac {n}{pd}}) \pmod {p^{\nu_p(n)}}.$
\end{lem}

\subsection{$q$-Deformation of Witt rings}
\label{definition of addition and multi}

Let $\textsf{R}=\mathbb Q[q][x_n,y_n: n\ge 1].$
Consider the map
$${\overline \Phi}^{q}_{\textsf R}: \textsf R^{\mathbb N} \to \textsf R^{\mathbb N},\quad
(a_n)_{n\in \mathbb N} \mapsto \left(\sum_{d|n}d (1-{q}^{\frac
{n}{d}})a_d^{\frac nd}\right)_{n\in \mathbb N}.$$
Let $X=(x_n)_{n\in \mathbb N}$ and $Y=(y_n)_{n\in \mathbb N}$
and then define $\textsf S_n, \textsf P_n \,\, (n\ge 1)$ via the following relations:
\begin{equation}\label{multiplication defining polynomial:witt ring}
\begin{aligned}
&\sum_{d|n}d \left[\frac nd\right]_{q}\textsf S_d^{\frac nd}=\left(\sum_{d|n}d \left[\frac nd \right]_{q}x_d^{\frac nd}\right)+
\left(\sum_{d|n}d\left[\frac nd\right]_{q}y_d^{\frac nd}\right),\\
&\sum_{d|n}d\left[\frac nd\right]_{q}\textsf P_d^{\frac nd}=(1-q)\left(\sum_{d|n}d \left[\frac nd\right]_{q}x_d^{\frac nd}\right)
\left(\sum_{d|n}d \left[\frac nd\right]_{q}y_d^{\frac nd}\right).
\end{aligned}
\end{equation}
Multiplying either side of Eq.\eqref{multiplication defining polynomial:witt ring} by $1-q$, one can see that
\begin{align*}
&{\overline \Phi}^{q}_{\textsf R}((\textsf S_n)_{n\in \mathbb N})={\overline \Phi}^{q}_{\textsf R}(X)
+{\overline \Phi}^{q}_{\textsf R}(Y), \text{ and }\\
&{\overline \Phi}^{q}_{\textsf R}((\textsf P_n)_{n\in \mathbb N})={\overline \Phi}^{q}_{\textsf R}(X)\cdot {\overline \Phi}^{q}_{\textsf R}(Y).
\end{align*}
We also let $\textsf I=(\textsf I_n)_{n\in \mathbb N}$ be the unique element of $\textsf R^{\mathbb N}$
satisfying the condition that
\begin{equation}\label{how to define additive inverse}
\begin{aligned}
&\left(\sum_{d|n}d \left[\frac nd\right]_{q}x_d^{\frac nd}\right)+
\left(\sum_{d|n}d\left[\frac nd\right]_{q}\textsf I_d^{\frac nd}\right)=0 \quad (n\ge 1).
\end{aligned}
\end{equation}
Let us multiply either side of Eq.\eqref{how to define additive inverse} by $1-q$ to get
$${\overline \Phi}^{q}_{\textsf R}(X)+ {\overline \Phi}^{q}_{\textsf R}(\textsf I)=0.$$
Utilizing the mathematical induction on $n$,
one can easily see that
$\textsf S_n,\textsf P_n \in \mathbb Q[q][x_d,y_d:d|n]$
and $\textsf I_n \in \mathbb Q[q][x_d:d|n]$ for all $n\ge 1.$
The following lemma is a key ingredient in the construction of ${\overline {W}}^{q}$.

\begin{lem}\label{key lemma in the proof of existence of q deformation of Witt rings}
For each positive integer $n$, it holds that
$\textsf S_n,\textsf P_n\in \mathbb Z[q][x_d,y_d:d|n]$ and
$\textsf I_n \in \mathbb Z[q][x_d:d|n].$
\end{lem}

\pf
We'll only prove that $\textsf P_n\in \mathbb Z[q][x_d,y_d:d|n]$
since the others can be treated in the same way.
When $n=1$, then our assertion is obvious since $P_1=(1-q)x_1y_1$.
If $n>1$, then assume that our assertion holds for all positive integers less than $n$.
Note that
\begin{equation}\label{multiplication defining polynomial:witt ring 1}
\begin{aligned}
&n(1-q)\textsf P_n=
\left(\sum_{d|n}d (1-q^{\frac
{n}{d}})x_d^{\frac nd}\right)
\left(\sum_{d|n}d (1-q^{\frac{n}{d}})y_d^{\frac nd}\right)
-\sum_{d|n \atop d<n}d (1-q^{\frac
{n}{d}})\textsf P_d^{\frac nd}, \quad (n\ge 1).
\end{aligned}
\end{equation}
First, we note that the induction hypothesis shows that the right hand side of Eq.\eqref{multiplication defining polynomial:witt ring 1}
is contained in $(1-q)\mathbb Z[q][x_d,y_d:d|n].$
Second, let us show that it is divisible by $n$.
We will do this by using the mathematical induction on $n$.
Let $p$ be any prime divisor of $n$.
Ignoring all the terms with $\nu_p(d)=\nu_p(n)$, one can derive the following modular equivalence:
\begin{equation}\label{reduced modular equivalence of new Witt rings}
\begin{aligned}
&n(1-q)\textsf P_n\\
&\equiv \left(\sum_{d|\frac np }d (1-q^{\frac
{n}{d}})x_d^{\frac {n}{d}}\right)
\left(\sum_{d|\frac np }d (1-q^{\frac{n}{d}})y_d^{\frac {n}{d}}\right)
-\sum_{d|\frac np}d (1-q^{\frac
{n}{d}})\textsf P_d^{\frac {n}{d}} \pmod {p^{\nu_p(n)}}\\
&=\left(\sum_{d|\frac np }d(1-(q^p)^{\frac
{n}{pd}})(x_d^p)^{\frac {n}{pd}}\right)
\left(\sum_{d|\frac np }d (1-(q^p)^{\frac{n}{pd}})(y_d^p)^{\frac {n}{pd}}\right)
-\sum_{d|\frac np}d (1-(q^p)^{\frac{n}{pd}})(\textsf P_d^p)^{\frac {n}{pd}}.
\end{aligned}
\end{equation}
In view of Eq.\eqref{multiplication defining polynomial:witt ring 1}, one can see that the term
$$\left(\sum_{d|\frac np } d (1-(q^p)^{\frac{n}{pd}})(x_d^p)^{\frac {n}{pd}}\right)
\left(\sum_{d|\frac np } d (1-(q^p)^{\frac{n}{pd}})(y_d^p)^{\frac {n}{pd}}\right)$$
is equal to
$$\sum_{d|\frac np}d (1-(q^p)^{\frac{n}{pd}})
(\textsf P_d|_{q \mapsto q^p \atop {x_i \mapsto x_i^p \atop y_i \mapsto y_i^p}})^{\frac {n}{pd}},$$
where $\textsf P_d|_{q \mapsto q^p \atop {x_i \mapsto x_i^p \atop y_i \mapsto y_i^p}}$ denotes the polynomial obtained from $\textsf P_d$
by replacing $q, x_i,y_i$ by $q^p,x_i^p,y_i^p$, respectively, for every divisor $i$ of $d$.
Therefore, Eq.\eqref{reduced modular equivalence of new Witt rings} can be simplified as
\begin{equation}\label{reduction of the defining relation}
n(1-q)\textsf P_n\equiv \sum_{d|\frac np}d (1-(q^p)^{\frac{n}{pd}})
\left((\textsf P_d|_{q \mapsto q^p \atop {x_i \mapsto x_i^p \atop y_i \mapsto y_i^p}})^{\frac {n}{pd}}
-(\textsf P_d^p)^{\frac {n}{pd}}\right) \pmod {p^{\nu_p(n)}}.
\end{equation}
For each divisor $d$ of $n/p$,
the induction hypothesis implies that $\textsf P_d \in \mathbb Z[q][x_i,y_i: i|d]$.
Hence, it follows from Lemma \ref{key modular equivalence 1}(a) that
$$(\textsf P_d|_{q \mapsto q^p \atop {x_i \mapsto x_i^p \atop y_i \mapsto y_i^p}})^{\frac {n}{pd}}
-\textsf P_d^{\frac {n}{d}}\equiv 0 \pmod{p^{\nu_p(n)-\nu_p(d)}}.$$
So we can conclude that the right hand side of Eq.\eqref{reduction of the defining relation} is divisible by $p^{\nu_p(n)}$.
Since $p$ is arbitrary, it follows that $n(1-q)\textsf P_n$ is divisible by $n$.
Because $1-q$ being a primitive polynomial, this implies that $\textsf P_n \in \mathbb Z[q][x_d,y_d:d|n]$.

In the same manner, we can show that
$\textsf S_n\in \mathbb Z[q][x_d,y_d:d|n].$
To be more precise,
$$\textsf S_n-(x_n+y_n) \in \mathbb Z[q][x_d,y_d:d|n, d<n], \quad (n\ge 1)$$
since the coefficient of $x_n,y_n$ in $\textsf  S_n$ equals 1, respectively.
Combining the identity
$$\textsf  S_n|_{x_d \mapsto x_d \atop y_d \mapsto \textsf I_d}=0$$
with the induction hypothesis that $\textsf  I_d \in \mathbb Z[q][x_i\,:\,i|d]$ for all positive integers $d$ less than $n$,
we finally conclude that
$\textsf I_n +x_n \in \mathbb Z[q][x_d :d|n \text{ with } d<n].$
This completes the proof.
\qed

\vskip 2mm
Given a polynomial $g(q)\in \mathbb Z[q]$, the general ${{{\overline {W}}}}^{g(q)}$ can be obtained
merely by transport of $\mathbb Z[q]$-algebra structure by the map $\mathbb Z[q] \to \mathbb Z[q]$ sending $q$ to $g(q)$.
Set $\textsf{gh}:=\iota^{\textsf{Rings}}_{\textsf{RINGS}}\circ gh.$
With this notation, we can derive the following generalization of Witt vectors.

\begin{thm}\label{Construction of new Witt rings attached to polynomials}
Given any polynomial $g(q)\in \mathbb Z[q]$, there exists a functor
$${{{\overline {W}}}}^{g(q)}:\textsf{$\mathbb Z[q]$-Algebras} \to \textsf{RINGS}$$
subject to the following conditions$:$
\begin{enumerate}
\item
As a set, ${{{\overline { W}}}}^{g(q)}(A)=A^\mathbb N$.
\item
For every $\mathbb Z[q]$-algebra homomorphism $f:A\to B$
and every $\alpha \in {{{\overline { W}}}}^{g(q)}(A)$ one has
${{{\overline { W}}}}^{g(q)}(f)(\alpha)=f\circ \alpha$.
\item
The map ${\overline \Phi}^{g(q)}_{A}:{{{\overline {W}}}}^{g(q)}(A)\to \textsf{gh}(A)$ defined by
$$(a_n)_{n\in \mathbb N} \mapsto \left(\sum_{d|n}d(1-g(q)^{\frac nd})a_d^{\frac nd}\right)_{n\in \mathbb N}$$
is a ring homomorphism.
\end{enumerate}
If $g(q)\ne 1$, then ${{{\overline {W}}}}^{g(q)}$ is uniquely determined by the conditions (1)-(3).
\end{thm}

\pf
First, let us show the existence of ${{\overline {W}}}^{g(q)}$.
Let $\textsf F=\mathbb Z[q][x_n,y_n: n\ge 1]$,
${X}=(x_n)_{n\in \mathbb N}$ and ${Y}=(y_n)_{n\in \mathbb N}$.
Given any $\mathbb Z[q]$-algebra $A$, define ${{\overline {W}}}^{g(q)}(A)$ to be the ring $A^{\mathbb N}$ whose ring operations are defined in
such a way that
\begin{align*}
&(a_n)_{n\in \mathbb N}+(b_n)_{n\in \mathbb N}= (\textsf S_n|_{x_d \mapsto a_d\atop y_d \mapsto b_d})_{n\in \mathbb N},\\
&(a_n)_{n\in \mathbb N}\cdot (b_n)_{n\in \mathbb N}= (\textsf P_n|_{x_d \mapsto a_d\atop y_d\mapsto b_d})_{n\in \mathbb N},\\
&-(a_n)_{n\in \mathbb N}=(\textsf I_n|_{x_d \mapsto a_d})_{n\in \mathbb N}
\end{align*}
for all elements $(a_n)_{n\in \mathbb N},(b_n)_{n\in \mathbb N} \in {{\overline {W}}}^{g(q)}(A)$.
In addition, define morphisms as in (2).
Then one easily sees that the assignment $A\mapsto {{\overline {W}}}^{g(q)}(A)$
induces a functor
$${{\overline {W}}}^{g(q)}: \textsf{$\mathbb Z[q]$-Algebras} \to \textsf{RINGS}$$
satisfying the conditions (1)-(3).

For the second assertion, assume that $g(q)\ne 1$ and $\hat {W}^{g(q)}$ is
a functor satisfying the desired conditions (1)-(3), where
the third condition says that the map
$${{\hat \Phi}}^{g(q)}_{A}:{\hat {W}}^{g(q)}(A)\to \textsf{gh}(A),\,\,
(a_n)_{n\in \mathbb N} \mapsto
\left(\sum_{d|n}d(1-g(q)^{\frac nd})a_d^{\frac nd}\right)_{n\in \mathbb N}$$
is a ring homomorphism.
For any $\mathbb Z[q]$-algebra $A$, let
${\rm id}_A:{{\overline {W}}}^{g(q)}(A)\to \hat {W}^{g(q)}(A)$ be the identity map.
We claim that $A \mapsto {\rm id}_A$ induces a natural isomorphism ${\rm id}:{{\overline {W}}}^{g(q)}\to \hat {W}^{g(q)}.$
In verifying our claim, the only nontrivial part is to confirm that ${\rm id}_A$ is a ring homomorphism.
In case where $A=\textsf F$, this follows from the fact that
${\overline \Phi}^{g(q)}_{\textsf F}, \hat {\overline \Phi}^{g(q)}_{\textsf F}$ are injective
and ${\overline \Phi}^{g(q)}_{\textsf F}={\hat \Phi}^{g(q)}_{\textsf F}\circ {\rm id}_{\textsf F} $.
From now on, suppose that $A$ is an arbitrary $\mathbb Z[q]$-algebra.
Given any elements ${a}=(a_n)_{n\in \mathbb N}$ and ${b}=(b_n)_{n\in \mathbb N}$ of ${{\overline {W}}}^{g(q)}(A)$, note that
$${\rm id}_A({a}\cdot {b})={\rm id}_A((\textsf P_n|_{x_d \mapsto a_d\atop y_d \mapsto b_d})_{n\in \mathbb N})
=(\textsf P_n|_{x_d \mapsto a_d\atop y_d \mapsto b_d})_{n\in \mathbb N}.$$
We now compute ${\rm id}_A({a})\cdot {\rm id}_A(b)$.
Let $\pi: \textsf F\to A$ be the unique $\mathbb Z[q]$-algebra homomorphism with $x_n\mapsto a_n, \,\,y_n \mapsto b_n$
for all $n\in \mathbb N$.
The functoriality of $\hat {W}^{g(q)}$ induces the ring homomorphism
$\hat {W}^{g(q)}(\pi):\hat {W}^{g(q)}(\textsf F) \to \hat {W}^{g(q)}(A)$, which satisfies that
$X=(x_n) \mapsto a$ and $Y=(y_n) \mapsto b$.
Therefore, it follows that
$$
{\rm id}_A({a})\cdot {\rm id}_A(b)
=a\cdot b
=\hat {W}^{g(q)}(\pi)({X})\cdot \hat {W}^{g(q)}(\pi)({Y})
=\hat {W}^{g(q)}(\pi)({X}\cdot {Y}).$$
Finally we replace ${X}\cdot {Y}$ by $(\textsf P_n)_{n\in \mathbb N}$ to derive that
$${\rm id}_A({a})\cdot {\rm id}_A(b)
=(\pi(\textsf P_n))_{n\in \mathbb N}
=(\textsf P_n|_{x_d \mapsto a_d\atop y_d \mapsto b_d})_{n\in \mathbb N}.$$
In the same way as in the above, we can verify that
$${\rm id}_A({a}+{b})=(\textsf S_n|_{x_d \mapsto a_d\atop y_d \mapsto b_d})_{n\in \mathbb N}
={\rm id}_A({a})+{\rm id}_A({b}).$$
This completes the proof.
\qed

\begin{rem}\label{shapening of the third condition}
{\rm

(a) One easily sees that the assignment $A\mapsto {\overline \Phi}^{g(q)}_A$ induces a natural transformation
${\Phi}^{g(q)}:{{{\overline {W}}}}^{g(q)}\to \textsf{gh}\circ F^{\textsf{$\mathbb Z[q]$-Algebras}}_{\textsf{Rings}}$.

(b) The third condition of Theorem \ref{Construction of new Witt rings attached to polynomials} can be sharpened a little more.
Indeed the first equation in \eqref{multiplication defining polynomial:witt ring}
implies that the third condition can be divided into two parts in the following manner:

\noindent
{\it (3-1)
The map $(\Phi_{\rm ab})^{g(q)}_{A}: {{\overline {W}}}^{g(q)}(A)\to \textsf{gh}(A)$ defined by
$$(a_n)_{n\in \mathbb N}\mapsto \left(\sum_{d|n}d\left[\frac nd \right]_{g(q)} a_d^{\frac nd}\right)_{n\in \mathbb N}$$
is an (abelian) group homomorphism.

\noindent
(3-2)
The map ${\overline \Phi}^{g(q)}_{A}: {{\overline {W}}}^{g(q)}(A)\to \textsf{gh}(A)$ defined by
$$(a_n)_{n\in \mathbb N} \mapsto
\left(\sum_{d|n}d(1-{g(q)}^{\frac nd})a_d^{\frac nd}\right)_{n\in \mathbb N}$$
is multiplicative.}\\
This description is very useful in studying the abelian structure of ${{\overline {W}}}^{g(q)}(A)$.
Furthermore, by mimicking the proof of Theorem \ref{Construction of new Witt rings attached to polynomials}, one can show that
$F^{\textsf{RINGS}}_{\textsf {Abel}}\circ {{{\overline {W}}}}^{g(q)}$ is the unique functor satisfying the conditions
(1), (2), and (3-1) for all $g(q)\in \mathbb Z[q]$.
In this case, it should be mentioned that the condition $g(q)\ne 1$ is not necessary in proving the uniqueness.
}\end{rem}

\subsection{Specializing $q$ into an integer}
It is quite interesting to note that one can specialize $q$ into an arbitrary integer $m$ in Theorem \ref{Construction of new Witt rings attached to polynomials}.
In fact this specialization makes sense since it is compatible with ring homomorphisms.
Let us explain it in more detail.
Let $A,B$ be $\mathbb Z[q]$-algebras and let $f: A\to B$ be a $\mathbb Z[q]$-algebra homomorphism.
Denote by $q\mapsto m, \,\,A|_{q\mapsto m},\,\, B|_{q\mapsto m}, \text{ and }f|_{q\mapsto m}$
the specialization of $q$ into $m$, the ring obtained from $A$ by the specialization $q\mapsto m$, the ring obtained from $B$ by the specialization $q\mapsto m$, and the ring homomorphism obtained from $f$ by the specialization $q\mapsto m$, respectively.
Then the following diagram
\begin{equation*}
\begin{CD}
A@>f>>B\\
@V{q \mapsto m}VV  @VV q \mapsto m V \\
A|_{q\mapsto m}@>>f|_{q\mapsto m}>B|_{q\mapsto m}
\end{CD}
\end{equation*}
commutes since
$$(({q \mapsto m}) \circ f)(q\cdot x)=({q \mapsto m})(q\cdot f(x))=mf(x)$$
and
$$( f_{q \mapsto m}\circ ({q \mapsto m}))(q\cdot x)= f_{q \mapsto m}(mx)=mf(x)$$
for all $x\in A$.
Denote by
$${W}^m:\textsf{Rings} \to \textsf{RINGS}$$
the functor obtained from ${{\overline {W}}}^q$
by specializing $q$ into $m$.
Let us explain it in more detail.
For any commutative ring $A$ with unity,
${W}^m(A)$ means the ring defined in such a way that
\begin{align*}
&(a_n)_{n\in \mathbb N}+(b_n)_{n\in \mathbb N}
= (\textsf S_n|_{{q\mapsto m \atop x_d \mapsto a_d}\atop y_d \mapsto b_d})_{n\in \mathbb N},\\
&(a_n)_{n\in \mathbb N}\cdot (b_n)_{n\in \mathbb N}
=(\textsf P_n|_{{q\mapsto m \atop x_d \mapsto a_d}\atop y_d \mapsto b_d})_{n\in \mathbb N},\\
&-(a_n)_{n\in \mathbb N}=(\textsf I_n|_{{q\mapsto m \atop x_d \mapsto a_d}})_{n\in \mathbb N}
\end{align*}
for any elements $(a_n)_{n\in \mathbb N},(b_n)_{n\in \mathbb N}\in {W}^m(A)$.
The following corollary follows immediately from Theorem \ref{Construction of new Witt rings attached to polynomials}.

\begin{cor}\label{Witt rings attached to integers}
Given an integer $m$, the functor ${W}^m$ satisfies
the following conditions$:$
\begin{enumerate}
\item
As a set, ${W}^m(A)=A^\mathbb N$.
\item
For every unity-preserving ring homomorphism $f:A\to B$
and every $\alpha \in {W}^m(A)$ one has
${W}^m(f)(\alpha)=f \circ \alpha$.
\item
The map $\Phi^m_{A}: {W}^m(A)\to \textsf{gh}(A)$ defined by
$$(a_n)_{n\in \mathbb N} \mapsto
\left(\sum_{d|n}d(1-m^{\frac nd})a_d^{\frac nd}\right)_{n\in \mathbb N}$$
is a ring homomorphism.
\end{enumerate}
Furthermore, if $m\ne 1$, then it is uniquely determined by (1)-(3).
\end{cor}

\begin{rem}{\rm

(a) As in Remark \ref{shapening of the third condition},
the third condition in
Corollary \ref{Witt rings attached to integers} can be rephrased in the following manner:

{\it (3-1)
The map $(\Phi_{\rm ab})^{m}_{A}: {W}^{m}(A)\to \textsf{gh}(A)$ defined by
$$(a_n)_{n\in \mathbb N} \mapsto
\left(\sum_{d|n}d\left[\frac nd \right]_{m} a_d^{\frac nd}\right)_{n\in \mathbb N}$$
is an (abelian) group homomorphism.

(3-2)
The map ${\overline \Phi}^{g(q)}_{A}: {W}^{m}(A)\to \textsf{gh}(A)$ defined by
$$(a_n)_{n\in \mathbb N} \mapsto
\left(\sum_{d|n}d(1-{g(q)}^{\frac nd})a_d^{\frac nd}\right)_{n\in \mathbb N}$$
is multiplicative.}

(b) Let $m$ be an integer.
In the same way as in Corollary \ref{Witt rings attached to integers}, we can derive a functor
$${\overline {W}}^{g(q)}|_{q\mapsto m}: \textsf{Rings} \to \textsf{RINGS}$$
by specializing $q$ into $m$.
Let us denote this functor by $W^{g(m)}$.
To see that this notation is well defined,
let us assume that $g(m)=g'(m')=k$
for some $g(q),g'(q)\in \mathbb Z[q]-\{1\}$ and
$m,m'\in \mathbb Z$.
Then it is obvious that
$$W^{k}(\stackrel{\text{by def.}}{=}{\overline {W}}^{q}|_{q\mapsto k})={\overline {W}}^{g(q)}|_{q\mapsto m}
={\overline {W}}^{g'(q)}|_{q\mapsto m'}.$$

(c) It should be noted that ${\overline {W}}^m$ is not $W^m$ but ${\overline {W}}^{g(q)}$ with $g(q)=m$.
}\end{rem}

\subsection{Classification theorems}
\label{Classification Theorem( Witt ring)}

\subsubsection{}
We classify ${\overline {W}}^{g(q)}$ and $F^{\textsf{RINGS}}_{\textsf {Abel}}\circ{{{\overline {W}}}}^{g(q)}$
up to strict-isomorphism as $g(q)$ varies over $\mathbb Z[q]$.
To begin with, we start with the definition of strict-isomorphism.

\begin{df} Suppose that $g(q)$ and $h(q)$ range over $\mathbb Z[q]$.

{\rm (a)} We say that ${{{\overline {W}}}}^{g(q)}$ is strictly isomorphic to ${{{\overline {W}}}}^{h(q)}$
if there exists a natural isomorphism $\omega:{{{\overline {W}}}}^{g(q)}\to {{{\overline {W}}}}^{h(q)}$ such that
${\overline \Phi}^{g(q)}={\overline \Phi}^{h(q)}\circ \omega$.

{\rm (b)} We say that $F^{\textsf{RINGS}}_{\textsf {Abel}}\circ{{{\overline {W}}}}^{g(q)}$ is strictly isomorphic to
$F^{\textsf{RINGS}}_{\textsf {Abel}}\circ{{{\overline {W}}}}^{h(q)}$
if there exists a natural isomorphism $\omega:F^{\textsf{RINGS}}_{\textsf {Abel}}\circ{{{\overline {W}}}}^{g(q)}\to F^{\textsf{RINGS}}_{\textsf {Abel}}\circ{{{\overline {W}}}}^{h(q)}$ such that
${\overline \Phi}^{g(q)}_{\rm ab}={\overline \Phi}^{h(q)}_{\rm ab}\circ \omega$.
\end{df}

\begin{prop}\label{classification of Witt rings attached to polynomials}
Let $g(q),h(q) \in \mathbb Z[q]$. Then
${{{\overline {W}}}}^{g(q)}$ is strictly isomorphic to ${{{\overline {W}}}}^{h(q)}$
if and only if $h(q)=g(q)$ or $h(q)=2-g(q)$ with $g(q)\ne 1$.
\end{prop}

\pf
Assume that ${{{\overline {W}}}}^{g(q)}$ is strictly isomorphic to ${{{\overline {W}}}}^{h(q)}$.
Let
$F=\mathbb Z[q][x_n: n\ge 1].$
Then there exists a ring isomorphism
$\omega_F:{{{\overline {W}}}}^{g(q)}(F)\to {{{\overline {W}}}}^{h(q)}(F) $
satisfying that
\begin{equation}\label{isom for classification}
\sum_{d|n}d(1-g(q)^{\frac nd})x_d^{\frac nd}=\sum_{d|n}d(1-h(q)^{\frac nd})\textsf{y}_d^{\frac nd}, \quad (n\ge 1),
\end{equation}
where $( \textsf{y}_n)_{n\in \mathbb N}:=\omega_F((x_n)_{n\in \mathbb N}).$
If $g(q)=1$, then
$$\sum_{d|n}d(1-h(q)^{\frac nd})\textsf{y}_d^{\frac nd}=0, \quad (n\ge 1).$$
If $h(q)\ne 1$, then $\textsf{y}_n=0$ for all $n\ge 1$. This is obviously a contradiction
since $\omega_F$ is a bijection.
Thus $h(q)$ should be $1$.
From now on, assume that $g(q)\ne 1$.
When $n=1$, Eq.\eqref{isom for classification} is reduced to
$$(1-g(q))x_1=(1-h(q))\textsf{y}_1.$$
Hence, $\textsf{y}_1$ is divisible by $x_1$.
If we let $y_1=x_1\alpha$ for some $\alpha\in F,$ then
$$(1-g(q))x_1=(1-h(q))x_1 \alpha$$
and thus $1-g(q)$ is divisible by $1-h(q).$
In the same manner, we can show that $1-h(q)$ is divisible by $1-g(q).$
As a consequence, $1-g(q)=\pm (1-h(q)).$
If $g(q)\ne h(q)$, then $1-g(q)=-(1-h(q))$, that is, $h(q)=2-g(q)$.
So we are done.

Next, let us prove the converse.
In case where $h(q)=g(q)$ we can take $\omega$ as ${\rm id}$.
From now on, we suppose $h(q)\ne g(q)$ .
To begin with, note $g(q)\ne 1$ since $h(q)=2-g(q)=1$.
Also we note $1-h(q)=-(1-g(q)).$
As before, let
$F=\mathbb Z[q][x_n: n\ge 1]$
and ${X}=(x_n)_{n\in \mathbb N}$.
For each positive integer $n$, define $\textsf{y}_n$
via the following recursive relations:
\begin{equation}\label{classification of Witt rings; defining polynomial}
-n(1-g(q))\textsf{y}_n=\sum_{d|n}d(1-g(q)^{\frac nd})x_d^{\frac nd}-\sum_{d|n \atop d<n}d(1-(2-g(q))^{\frac nd})\textsf{y}_d^{\frac nd},\quad (n\ge 1).
\end{equation}
We claim that
$$\textsf{y}_n\in \mathbb Z[g(q)][x_d: d|n], \quad (n\ge 1).$$
When $n=1$, it is clear that $\textsf y_1=-x_1$.
If $n>1$, assume that our assertion holds for all positive integers less than $n$.
Exploiting the induction hypothesis, it is easy to see that the right hand side of Eq.\eqref{classification of Witt rings; defining polynomial}
is contained in $(1-g(q))\mathbb Z[g(q)][x_d: d|n]$.
In the following, we will see that it is also divisible by $n$.
Choose any prime divisor $p$ of $n$.
Note that
\begin{align*}
-n(1-g(q))\textsf{y}_n
\equiv \sum_{d| \frac np}d(1-(g(q)^p)^{\frac {n}{pd}})(x_d^p)^{\frac {n}{pd}}
-\sum_{d| \frac np}d(1-(2-g(q))^{\frac {n}{d}})(\textsf{y}_d^p)^{\frac {n}{pd}} \pmod {p^{\nu_p(n)}}.
\end{align*}
But, by replacing $n$ and $g(q)$ by $n/p$ and $g(q)^p$, respectively, in Eq.\eqref{classification of Witt rings; defining polynomial},
we can deduce that
$$\sum_{d| \frac np}d(1-(g(q)^p)^{\frac {n}{pd}})(x_d^p)^{\frac {n}{pd}}
=\sum_{d| \frac np}d(1-(2-g(q)^p)^{\frac {n}{pd}})(\textsf{y}_d|_{g(q) \mapsto g(q)^p \atop x_i \mapsto x_i^p})^{\frac {n}{pd}}.$$
Consequently we have
\begin{align*}
&-n(1-g(q))\textsf{y}_n\\
&\equiv \sum_{d| \frac np}d(1-(2-g(q)^p)^{\frac {n}{pd}})(\textsf{y}_d|_{g(q) \mapsto g(q)^p \atop x_i \mapsto x_i^p})^{\frac {n}{pd}}
-\sum_{d| \frac np}d(1-(2-g(q))^{\frac {n}{d}})(\textsf{y}_d^p)^{\frac {n}{pd}} \pmod {p^{\nu_p(n)}}\\
&\equiv \sum_{d| \frac np}d(1-(2-g(q)^p)^{\frac {n}{pd}})((\textsf{y}_d|_{g(q) \mapsto g(q)^p \atop x_i \mapsto x_i^p})^{\frac {n}{pd}}-(\textsf{y}_d^p)^{\frac {n}{pd}} ) \\
&\,\,-\sum_{d| \frac np}d((2-g(q)^p)^{\frac {n}{pd}}-((2-g(q))^p)^{\frac {n}{pd}})(\textsf{y}_d^p)^{\frac {n}{pd}} \pmod {p^{\nu_p(n)}}.
\end{align*}
For each divisor $d$ of $n/p$, Lemma \ref{key modular equivalence 1}(a) implies that
$$(\textsf{y}_d|_{g(q) \mapsto g(q)^p \atop x_i \mapsto x_i^p})^{\frac {n}{pd}}-(\textsf{y}_d^p)^{\frac {n}{pd}} \equiv 0 \pmod {p^{\nu_p(n)-\nu_p(d)}}$$
and
$$(2-g(q)^p)^{\frac {n}{pd}}-((2-g(q))^p)^{\frac {n}{pd}} \equiv 0 \pmod {p^{\nu_p(n)-\nu_p(d)}}$$
since
$$(2-g(q)^p)-(2-g(q))^p\equiv 0 \pmod p.$$
Therefore, it holds that
$$-n(1-g(q))\textsf{y}_n\equiv 0 \pmod {p^{\nu_p(n)}}$$
and thus $-n(1-g(q))\textsf{y}_n\equiv 0 \pmod {n}$
since $p$ is an arbitrary prime divisor of $n$.
For the full justification of our claim,
we should prove that the right hand side of Eq.\eqref{classification of Witt rings; defining polynomial} is divisible by
$n(1-g(q))$.
This can be done in the same way as in the proof of Lemma \ref{key lemma in the proof of existence of q deformation of Witt rings}.

Now, we are ready to construct the desired natural transformation.
Given a $\mathbb Z[q]$-algebra $A$, consider the map
$$\omega_A:{{{\overline {W}}}}^{g(q)}(A)\to {{{\overline {W}}}}^{2-g(q)}(A), \,\,(a_n)_{n\in \mathbb N}\mapsto (\textsf{y}_n|_{x_d \mapsto a_d}).$$
Obviously $\omega_A$ is a bijection since the coefficient of $x_n$ in $\textsf y_n$ equals 1.
One can verify in a routine way that the assignment $A \mapsto \omega_A$ induces a natural isomorphism $\omega:{{{\overline {W}}}}^{g(q)}\to {{{\overline {W}}}}^{2-g(q)}$
satisfying
${\overline \Phi}^{g(q)}={\overline \Phi}^{h(q)}\circ \omega$.
This completes the proof.
\qed

\begin{lem}\label{lemma for classification abelian Witt}
Let $p$ be a prime. For any polynomial $f(q)\in \mathbb Z[q]$,
$p|1-f(q)$ if and only if $p|[p]_{f(q)}$.
\end{lem}

\pf
Suppose $p|1-f(q)$.
If $f(q)=1$, there is nothing to prove since $[p]_{f(q)}=p$.
If $f(q)\ne 1$, let us write $1-f(q)$ as $-p^a\hat f(q)$ with
$(p,\hat f(q))=1$ and $a\ge 1$.
Then the desired result follows from
the fact that
$$1-f(q)^p=1-(1+p^a\hat f(q))^p\equiv 0 \pmod{p^{a+1}}.$$
Conversely, assume $p|[p]_{f(q)}$.
Then $p|1-f(q)^p$ and thus $p|1-f(q^p)$.
Letting
$$f(q)=a_0+a_1q+\cdots +a_mq^m,$$
the above condition implies that
$1-a_0, a_1,a_2, \cdots, a_m$ are divisible by $p$. Consequently
$p$ divides $1-f(q)$.
\qed

\vskip 2mm
Recall that a nonzero polynomial $k(q)\in \mathbb Z[q]$ is called {\it primitive} if
the only integers that divide all coefficients of $k(q)$ at once are $\pm 1$.
With this definition, we can derive the following classification theorem.

\begin{prop}\label{classification of Witt rings attached to polynomials; abelian version}
Let $g(q),h(q) \in \mathbb Z[q]$. Then
$F^{\textsf{RINGS}}_{\textsf {Abel}}\circ{{{\overline {W}}}}^{g(q)}$ is strictly isomorphic to
$F^{\textsf{RINGS}}_{\textsf {Abel}}\circ{{{\overline {W}}}}^{h(q)}$
if and only if $1-g(q)=c_1k(q)$ and $1-h(q)=c_2k(q)$ for some $c_1,c_2 \in \mathbb Z$ and $k(q)\in \mathbb Z[q]$ such that
\begin{enumerate}
\item $c_1$ and $c_2$ have the same set of prime divisors, and
\item $k(q)$ is primitive.
\end{enumerate}
\end{prop}

\pf
To begin with, assume that $F^{\textsf{RINGS}}_{\textsf {Abel}}\circ{{{\overline {W}}}}^{g(q)}$ is strictly isomorphic to
$F^{\textsf{RINGS}}_{\textsf {Abel}}\circ{{{\overline {W}}}}^{h(q)}$.
Let
$F=\mathbb Z[q][x_n: n\ge 1].$
By assumption there exists a ring isomorphism
$\omega_F:{{{\overline {W}}}}^{g(q)}(F)\to {{{\overline {W}}}}^{h(q)}(F) $
satisfying that
$$
\sum_{d|n}d\left[\frac nd \right]_{g(q)}x_d^{\frac nd}=\sum_{d|n}d\left[\frac nd \right]_{h(q)}\textsf{y}_d^{\frac nd}, \quad (n\ge 1),
$$
where $(\textsf{y}_n)_{n\in \mathbb N}:=\omega_F((x_n)_{n\in \mathbb N}).$
For a prime $p$, we have
\begin{equation}\label{divisibility of p g(q) and h(q)}
\textsf{y}_p=x_p+\frac{[p]_{g(q)}-[p]_{h(q)}}{p}x_1^p.
\end{equation}
It implies that $[p]_{g(q)}-[p]_{h(q)}$ is divisible by $p$.
Set
$$1-g(q)=c_1 q^a {\hat g}(q)\text{ and }1-h(q)=c_2 q^b {\hat h}(q),$$
where $c_1,c_2 \in \mathbb Z$, $a,b \in \mathbb Z_{\ge 0}$, and ${\hat g}(q), {\hat h}(q)$
are primitive polynomials with a positive constant term.
Then the assertion that $c_1,c_2$ have the same set of divisors follows from the following equivalences:
$$p|1-g(q) \stackrel{\text{by Lemma \ref{lemma for classification abelian Witt}}}{\Longleftrightarrow}
p|[p]_{g(q)} \stackrel{\text{by Eq.\ref{divisibility of p g(q) and h(q)}}}{\Longleftrightarrow}
p|[p]_{h(q)} \stackrel{\text{by Lemma \ref{lemma for classification abelian Witt}}}{\Longleftrightarrow} p|1-h(q)$$
Next, we will show that $a=b$ and $\hat g(q)=\pm \hat h(q)$.
In the following, suppose $p\nmid 1-g(q)$ (so, $p\nmid 1-h(q)$).
In particular, $c_1$ and $c_2$ are nonzero.
Because $g(q)^{p-1}-[p]_{h(q)}$ is divisible by $p$, so is
$$(1-h(q))g(q)^p-(1-g(q))(1-h(q)^p).$$
Hence, the congruences $g(q)^p\equiv g(q^p) \pmod p$ and $h(q)^p\equiv h(q^p) \pmod p$ imply that
$$(c_2q^b \hat h(q))(c_1q^{pa} \hat g(q^p))-(c_1q^a \hat g(q))(c_2q^{pb} \hat h(q^p))$$
is divisible by $p$.
For simplicity, set
$$\hat g(q)=a_0+a_1q+\cdots + a_sq^s$$
and
$$\hat h(q)=b_0+b_1q+\cdots + b_tq^t,$$
where $a_0, b_0$ are nonzero.
It is not difficult to see that
$$(c_2q^b \hat h(q))(c_1q^{pa} \hat g(q^p))-(c_1q^a \hat g(q))(c_2q^{pb} \hat h(q^p))$$
can be simplified to
$$c_1c_2q^{pa+b}\sum_{k\ge 0 }\left(\sum_{pi+j=k}a_ib_j \right)q^k-c_1c_2q^{a+pb}\sum_{k\ge 0}\left(\sum_{i+pj=k}a_ib_j \right)q^k.
$$
To see $a=b$, assume $a>b$.
Then the coefficient of $q^{a+pb}$ equals $c_1c_2a_0b_0$.
Since $p$ is an arbitrary prime with $p\nmid 1-g(q)$, we may assume that $p$ is sufficiently large.
Under this assumption, the condition $p|c_1c_2a_0b_0$ forces $a=0$ or $b=0$, which is absurd.
In the same manner, we can derive a contradiction when $a<b$.
Consequently $a=b$, thus
the coefficient of $q^{pa+a+k}$ equals $c_1c_2(a_0b_k-a_kb_0)$ for all $k\ge 0$
(because $p$ is sufficiently large).
In addition, since none of $c_1,c_2$ are divisible by $p$, we get $a_0b_k-a_kb_0=0$.
It means that $\hat h(q)=\frac {b_0}{a_0}\hat g(q)$.
On the other hand, $\hat g(q),\hat h(q)$ are primitive, thus $\frac {b_0}{a_0}$ should be $\pm 1$.
So we are done.

Next, let us prove the converse.
Let
$F=\mathbb Z[q][x_n: n\ge 1]$
and ${X}=(x_n)_{n\in \mathbb N}$.
For $n\ge 1$, define $\textsf{y}_n$
via the following recursive relations:
\begin{equation}\label{classification of abelianized Witt rings; defining polynomial}
n\textsf{y}_n=\sum_{d|n}d\left[\frac nd\right]_{g(q)}x_d^{\frac nd}-\sum_{d|n \atop d<n}d\left[\frac nd\right]_{h(q)}\textsf{y}_d^{\frac nd}
\quad (\forall n\ge 1).
\end{equation}
We claim that
$\textsf{y}_n\in \mathbb Z[g(q)][x_d: d|n]$ for all $n\ge 1.$
When $n=1$, $\textsf y_1=x_1$.
If $n>1$, assume that our assertion holds for all positive integers less than $n$.
In view of the induction hypothesis, it is obvious that the right hand side of Eq.\eqref{classification of abelianized Witt rings; defining polynomial}
is contained in $\mathbb Z[g(q)][x_d: d|n]$.
To see that it is also divisible by $n$,
choose any prime divisor $p$ of $n$.

First, assume that
$(p,1-g(q))=(p,1-h(q))=p.$
Using Lemma \ref{key modular equivalence 1}(a), one can easily see that
$d\left[\frac nd\right]_{g(q)}$ and $d\left[\frac nd\right]_{h(q)}$ are divisible by $p^{\nu_p(n)}$
for every divisor $d$ of $n$. As a result,
the right hand side of Eq.\eqref{classification of abelianized Witt rings; defining polynomial} is divisible by $p^{\nu_p(n)}$.

Second, assume that $(p,1-g(q))=(p,1-h(q))=1$.
From Eq.\eqref{classification of abelianized Witt rings; defining polynomial} it follows that
\begin{align*}
&n(1-g(q))(1-h(q))\textsf{y}_n\\
&\equiv (1-h(q))\sum_{d| \frac np}d(1-g(q)^{\frac nd})x_d^{\frac nd}-(1-g(q))\sum_{d| \frac np}d(1-h(q)^{\frac nd})\textsf{y}_d^{\frac nd}\\
&\equiv (1-h(q))\sum_{d| \frac np}d(1-(g(q)^p)^{\frac {n}{dp}})(x_d^p)^{\frac {n}{dp}}
-(1-g(q))\sum_{d| \frac np}d(1-(h(q)^p)^{\frac {n}{dp}})(\textsf{y}_d^p)^{\frac {n}{dp}}\\
&\equiv (1-h(q))(1-g(q^p))\sum_{d| \frac np}d\left[\frac {n}{dp}\right]_{g(q^p)}(x_d^p)^{\frac {n}{dp}}
-(1-g(q))(1-h(q^p))\sum_{d| \frac np}d\left[\frac {n}{dp}\right]_{h(q^p)}(\textsf{y}_d^p)^{\frac {n}{dp}}.
\end{align*}
Here `$\equiv$' denotes `$\equiv \pmod {p^{\nu_p(n)}}$'.
Considering Eq.\eqref{classification of abelianized Witt rings; defining polynomial} for $\frac np$, we can deduce that
$$\sum_{d| \frac np}d\left[\frac {n}{dp}\right]_{g(q^p)}(x_d^p)^{\frac {n}{dp}}=\sum_{d| \frac np}d\left[\frac {n}{dp}\right]_{h(q^p)}(\textsf{y}_d|_{q\mapsto q^p \atop x_d \mapsto x_d^p})^{\frac {n}{dp}}.$$
Combining it with the congruence
$$(\textsf{y}_d|_{q\mapsto q^p \atop x_d \mapsto x_d^p})^{\frac {n}{dp}}\equiv (\textsf{y}_d|_{q\mapsto q^p \atop x_d \mapsto x_d^p})^{\frac {n}{dp}}
\pmod{p^{\nu_p(\frac np)}},$$
we can see that
$n(1-g(q))(1-h(q))\textsf{y}_n$
is congruent to
$$((1-h(q))(1-g(q^p))-(1-g(q))(1-h(q^p)))
\sum_{d| \frac np}d\left[\frac {n}{dp}\right]_{h(q^p)}\textsf{y}_d^{\frac nd}
\pmod {p^{\nu_p(n)}},$$
which is zero.
Because $(p,1-g(q))=(p,1-h(q))=1$, it follows that
$$n\textsf{y}_n
\equiv 0 \pmod {p^{\nu_p(n)}}.$$
So our claim is verified.
Finally, by mimicking the method in the proof of Proposition \ref{classification of Witt rings attached to polynomials},
we can construct the desired natural transformation.
\qed

\subsubsection{}
Here we classify $W^{m}$
and $F^{\textsf{RINGS}}_{\textsf {Abel}}\circ W^{m}$ up to strict-isomorphism
as $m$ varies over the set of integers.

\begin{df}
Suppose that $a,b$ range over the set of integers.

{\rm (a)}
We say that $W^{a}$ is {\it strictly-isomorphic} to $W^{b}$
if there exists a natural isomorphism
$\omega: W^{a}\to W^{b}$ such that $\Phi^a=\Phi^b\circ \omega$.
In this case, $\omega$ is called a {\it strict natural isomorphism}.

{\rm (b)}
We say that $F^{\textsf{RINGS}}_{\textsf {Abel}}\circ W^{a}$
is strictly isomorphic to $F^{\textsf{RINGS}}_{\textsf {Abel}}\circ W^{b}$
if there exists a natural isomorphism
$\omega: F^{\textsf{RINGS}}_{\textsf {Abel}}\circ W^{a} \to F^{\textsf{RINGS}}_{\textsf {Abel}}\circ W^{b}$
such that $(\Phi_{\rm ab})^{a}=(\Phi_{\rm ab})^{b}\circ \omega$.
\end{df}

\begin{prop}\label{classifi of specialed witt ring}
{\rm (a)}
Let $a,b$ be integers.
Then $W^a$ is {\it strictly-isomorphic} to $W^b$
if and only if $a=b$ or $a+b=2$ with $a\ne 1$.

{\rm (b)}
Suppose that $m$ varies over the set of integers.
Then $W^m$ is classified by the set of positive integers up to strict-isomorphism.
\end{prop}

\pf
(a) It follows from proposition \ref{classification of Witt rings attached to polynomials}.

(b)
First, we will show that $W^a$ is not strictly isomorphic to $W^b$ for any integers $a,b \ge 1$ with $a\ne b$.
It is trivial by (a) since $a+b$ always exceeds 2.
Second, we will show that for any nonnegative integer $a$, $W^a$ is strictly isomorphic to $W^{b}$
for some positive integer $b$.
It, however, is trivial since
$W^a$ is strictly isomorphic to $W^{2-a}$ and $2-a\ge 2.$
\qed

\begin{prop}\label{abelian structure }
{\rm (a)}
Let $a,b$ be arbitrary integers.
Then $F^{\textsf{RINGS}}_{\textsf {Abel}}\circ W^a$
is strictly isomorphic to $F^{\textsf{RINGS}}_{\textsf {Abel}}\circ W^b$
if and only if $(1-a)$ and $(1-b)$ have the same set of prime divisors.

{\rm (b)}
Suppose that $m$ varies over the set of integers and let $\mathcal P$ be the set of all primes.
Then $F^{\textsf{RINGS}}_{\textsf {Abel}}\circ W^m$ is classified by $\{S \subset \mathcal P: |S|<\infty \text{ or } S=\mathcal P\}$ up to strict-isomorphism.
\end{prop}

\pf
(a) It follows from Proposition \ref{classification of Witt rings attached to polynomials; abelian version}.

(b) For $a,b\in \mathbb Z$ define $a\sim b$ if $1-a$ and $1-b$ have the same set of prime divisors.
Obviously $\sim$ is an equivalence relation on $\mathbb Z$.
And, from (a) it follows that $F^{\textsf{RINGS}}_{\textsf {Abel}}\circ W^m$ is classified by
the set of equivalence classes up to strict-isomorphism.
Consider the map
$$\eta: \mathbb Z/\sim \,\to\,  \{S \subset \mathcal P: |S|<\infty \text{ or } S=\mathcal P\}, \quad [m]\mapsto \text{ the set of all prime divisors of }1-m.$$
By the definition of $\sim$, $\eta$ is well defined and injective.
To see that $\eta$ is surjective, pick up any finite subset $\{p_1, \ldots, p_l\}$ of $\mathcal P$.
Then $\eta(1-p_1p_2\cdots p_l)=\{p_1, \ldots, p_l\}.$
In addition, $\eta(1)=\mathcal P$.
This completes the proof.
\qed

\section{$q$-Deformation of Grothendieck's $\Lambda$-functor}
\label{q-Deformation of Grothendieck's Lambda functor}

\subsection{$q$-Deformation of Grothendieck's $\Lambda$-functor and its specializations}
Let $g(q)$ be an arbitrary polynomial in $\mathbb Z[q]$.
For any $\mathbb Z[q]$-algebra $A$,
let ${\overline \Lambda}^{g(q)}(A)$ be the ring such that the underlying set consists of formal power series with coefficients in $A$ and constant 1.
For each object $A$ in \textsf{$\mathbb Z[q]$-Algebras}, the addition and the multiplication are defined in such a way that
\begin{equation}\label{How to define addition and multiplication g(q) case}
\begin{aligned}
&\prod_{n\ge 1}\frac{1}{1-a_nt^n} +_{{\overline \Lambda}^{g(q)}_A} \prod_{n\ge 1}\frac{1}{1-b_nt^n}
=\prod_{n\ge 1}\frac{1}{1-\textsf S_n|_{\atop {x_d\mapsto a_d \atop y_d \mapsto b_d}}t^n},\\
&\prod_{n\ge 1}\frac{1}{1-a_nt^n}\times_{{\overline \Lambda}^{g(q)}_A} \prod_{n\ge 1}\frac{1}{1-b_nt^n}
=\prod_{n\ge 1}\frac{1}{1-\textsf P_n|_{\atop {x_d\mapsto a_d \atop y_d \mapsto b_d}}t^n}.
\end{aligned}
\end{equation}
For the definition of $\textsf S_n, \textsf P_n$, see  Section \ref{definition of addition and multi}.
For a $\mathbb Z[q]$-algebra homomorphism $f:A\to B$, consider the map
${\overline \Lambda}^{g(q)}(f):{\overline \Lambda}^{g(q)}(A) \to {\overline \Lambda}^{g(q)}(B)$ defined by
$$ 1+\sum_{n\ge 1}a_nt^n \mapsto 1+\sum_{n\ge 1}f(a_n)t^n.$$
One easily sees that the assignment $A\mapsto {\overline \Lambda}^{g(q)}(A)$ induces a functor
$${\overline \Lambda}^{g(q)}: \textsf{$\mathbb Z[q]$-Algebras} \to \textsf{RINGS}.$$
Set $\textsf{Gh}:=\iota^{\textsf{Rings}}_{\textsf{RINGS}}\circ Gh$
(for the definition of $Gh$, refer to Section \ref{Witt ring functor and necklace ring functor}).
In the same way as in Theorem \ref{Construction of new Witt rings attached to polynomials}, 
we can derive the following characterization of ${\overline \Lambda}^{g(q)}$.

\begin{thm}\label{charcterization of q-deformed Grotendieck functor g(q)}
Let $g(q)$ be any polynomial in $\mathbb Z[q]$.
The functor ${\overline \Lambda}^{g(q)}$ satisfies the following conditions$:$
\begin{enumerate}
\item
As a set, ${\overline \Lambda}^{g(q)}(A)=\{1+\sum_{n\ge 1}a_nt^n : a_n \in A\}$.
\item
For a $\mathbb Z[q]$-algebra homomorphism $f:A\to B$, one has
$${\overline \Lambda}^{g(q)}(f):{\overline \Lambda}^{g(q)}(A) \to {\overline \Lambda}^{g(q)}(B),\quad 1+\sum_{n\ge 1}a_nt^n \mapsto 1+\sum_{n\ge 1}f(a_n)t^n.$$
\item
The map ${\overline \Upsilon}_A^{g(q)}: {\overline \Lambda}^{g(q)}(A) \to  \textsf{Gh}(A)$ defined by
$$\prod_{n\ge 1} \frac {1}{1-a_nt^n}
\mapsto \frac{d}{dt}\log \prod_{n\ge 1} \frac {1-g(q)a_nt^n}{1-a_nt^n}$$
is a ring homomorphism.
\end{enumerate}
Furthermore, if $g(q)\ne 1$, then ${\overline \Lambda}^{g(q)}$ is uniquely determined by the conditions (1)-(3).
\end{thm}

Consider the case where $g(q)=q$.
Specializing $q$ into an integer $m$ produces a functor
${\it \Lambda^m}: \textsf{Rings} \to \textsf{RINGS}$ for each integer $m.$
The subsequent corollary follows immediately from Theorem \ref{charcterization of q-deformed Grotendieck functor g(q)}.

\begin{cor}
For any integer $m$, $\Lambda^m$
satisfies the following conditions$:$
\begin{enumerate}
\item
As a set, $\Lambda^m(A)=\{1+\sum_{n\ge 1}a_nt^n : a_n \in A\}$.
\item
For a unity-preserving ring homomorphism $f:A\to B$, one has
$$\Lambda^m(f):\Lambda^m(A) \to \Lambda^m(B),\quad 1+\sum_{n\ge 1}a_nt^n \mapsto 1+\sum_{n\ge 1}f(a_n)t^n.$$
\item
The map $\Upsilon_A^m: \Lambda^m(A) \to \textsf{Gh}(A)$ defined by
$$\prod_{n\ge 1} \frac {1}{1-a_nt^n} \mapsto \frac{d}{dt}\log \prod_{n\ge 1} \frac {1-ma_nt^n}{1-a_nt^n}$$
is a ring homomorphism.
\end{enumerate}
Furthermore, if $m\ne 1$, then it is uniquely determined by (1)-(3).
\end{cor}

\subsection{Symmetric functions arising in the context of our $q$-deformed Witt rings.}
In this subsection, we focus on certain symmetric functions which are naturally arise in the context of our $q$-deformed Witt rings.
Let $X=\{x_i:i\ge 1 \}$ be an alphabet.
Define new symmetric functions $\textsf{u}_n^q(X),\textsf{v}_n^q(X), (n\ge 1)$ via the following relations:
\begin{equation}\label{rel between u and p}
\begin{aligned}
\sum_{d|n}d\left[\frac nd \right]_q\textsf{u}_d^q(X)^{\frac nd}=\sum_{d|n}d(1-q^{\frac nd})\textsf{v}_d^q(X)^{\frac nd}=p_n(X),\quad (n\ge 1).\\
\end{aligned}
\end{equation}
Here $p_n(X)$ denotes the $n$th power-sum symmetric function in $x_i$'s, that is,
$p_n(x)=x_1^n+x_2^n+x_3^n+\cdots.$
Let $H(X,t)$ be the generating function of $h_n(X)$'s,
that is,
$$H(X,t)=\sum_{n=0}^\infty h_n(X)t^n:=\prod_{i=1}^{\infty} \frac{1}{1-x_i t}.$$
Let
\begin{align*}
&\textsf{H}^q(X,t)=1+\sum_{n\ge 1}
\textsf{h}^q_n(X)t^n:=\prod_{n\ge 1}\frac{1}{1-\textsf{u}_n^q(X)t^n},  \text{ and }\\
&\textsf{G}^q(X,t)=1+\sum_{n\ge 1}\textsf{g}^q_n(X)t^n:=\prod_{n\ge 1}\frac{1}{1-\textsf{v}_n^q(X)t^n}.
\end{align*}

\begin{lem}
Let $X=\{x_i:i\ge 1 \},\, Y=\{y_i:i\ge 1 \}$ be alphabets.
For each positive integer $n$, we have

$\text{\rm (a) } \textsf{h}^q_n(X+Y)\in \mathbb Z[q][\textsf{h}^q_i(X),\textsf{h}^q_j(Y): 1 \le i,j \le n].$

$\text{\rm (b) } \textsf{g}^q_n(XY)\in \mathbb Z[q][\textsf{g}^q_i(X),\textsf{g}^q_j(Y): 1 \le i,j \le n].$
\end{lem}
\pf

(a)
Let us observe that
$$\textsf {u}_n^q(X)-\textsf{h}^q_n(X)\in \mathbb Z[q][\textsf{h}^q_d(X):d|n]$$
and
\begin{equation*}\label{relation between u and h}
\textsf{h}^q_n(X+Y)=\sum_{a_1+2a_2+\cdots+na_n=n \atop a_i \ge 0}
{\textsf {u}_1^q(X+Y)}^{a_1}{\textsf {u}_2^q(X+Y)}^{a_2}\cdots {\textsf {u}_n^q(X+Y)}^{a_n}
\end{equation*}
for all $n\ge 1$.
Thus, to prove our assertion,
it suffices to see that
$$\textsf {u}_n^q(X+Y) \in \mathbb Z[q][\textsf {u}_d^q(X),\textsf {u}_d^q(Y): d|n]$$ for all $n\ge 1.$
On the other hand, due to the identity
\begin{align*}
\sum_{d|n}d\left[\frac nd \right]_q\textsf {u}_d^q(X+Y)^{\frac nd}=p_n(X+Y)=p_n(X)+p_n(Y),
\end{align*}
we have
\begin{align*}
\sum_{d|n}d\left[\frac nd \right]_q\textsf {u}_d^q(X+Y)^{\frac nd}
=\sum_{d|n}d\left[\frac nd \right]_q(\textsf {u}_d^q(X)^{\frac nd}+\textsf {u}_d^q(Y)^{\frac nd}).
\end{align*}
Immediately it follows that
$$\textsf {u}_n^q(X+Y) \in \mathbb Q[q][\textsf {u}_d^q(X),\textsf {u}_d^q(Y): d|n]$$
for each positive integer $n$.
In the following, we will show that
$$\textsf {u}_n^q(X+Y) \in \mathbb Z[q][\textsf {u}_d^q(X),\textsf {u}_d^q(Y): d|n].$$
Let $p$ be any prime divisor of $n$.
Then
\begin{align*}
&n(1-q)\textsf {u}_n^q(X+Y)\\
&=\sum_{d|n}d(1-q^{\frac nd})(\textsf {u}_d^q(X)^{\frac nd}+\textsf {u}_d^q(Y)^{\frac nd})
-\sum_{d|n \atop d<n}d(1-q^{\frac nd})\textsf {u}_d^q(X+Y)^{\frac nd}\\
&\equiv \sum_{d| \frac np}d(1-(q^p)^{\frac {n}{pd}})
((\textsf {u}_d^q(X)^p)^{\frac {n}{dp}}+(\textsf {u}_d^q(Y)^p)^{\frac {n}{pd}})
-\sum_{d| \frac np}d(1-(q^p)^{\frac {n}{pd}})(\textsf {u}_d^q(X+Y)^p)^{\frac {n}{pd}}\\
&\equiv \sum_{d| \frac np}d(1-(q^p)^{\frac {n}{pd}})
((\textsf {u}_d^q(X+Y)|_{{q \mapsto q^p \atop \textsf {u}_i(X)\mapsto \textsf {u}_i(X)^p} \atop \textsf {u}_i(Y)\mapsto \textsf {u}_i(Y)^p})^{\frac {n}{pd}}-(\textsf {u}_d^q(X+Y)^p)^{\frac {n}{pd}})\\
&\equiv 0 \quad \text{(by Lemma \ref{key modular equivalence 1}).}
\end{align*}
Here `$\equiv$' denotes `$\equiv \pmod {p^{\nu_p(n)}}$'.
It implies that $n(1-q)\textsf {u}_n^q(X+Y)$ is divisible by $n$.
But, it it obvious that $n(1-q)\textsf {u}_n^q(X+Y)$ is divisible by $1-q$.
Hence, our proof is completed.

(b)
Recall
\begin{align*}
\textsf{g}^q_n(XY)=\sum_{a_1+2a_2+\cdots+na_n=n \atop a_i \ge 0}
{\textsf{v}_1^q(XY)}^{a_1}{\textsf{v}_2^q(XY)}^{a_2}\cdots {\textsf{v}_n^q(XY)}^{a_n}.
\end{align*}
Hence, to verify our assertion, it suffices to show that
$$\textsf{v}_n^q(XY) \in \mathbb Z[q][\textsf{v}_d^q(X),\textsf{v}_d^q(Y): d|n]$$
for all $n\ge 1.$
since
$$\textsf{v}_d^q(X)-\textsf{g}^q_d(X)\in \mathbb Z[q][g^q_i(X):i|d].$$
To this end, observe that the equalities
\begin{align*}
\sum_{d|n}d(1-q^{\frac nd})\textsf{v}_d^q(XY)^{\frac nd}=p_n(XY)=p_n(X)p_n(Y),\,\, (n\ge 1)
\end{align*}
imply
\begin{align*}
\sum_{d|n}d(1-q^{\frac nd})\textsf{v}_d^q(XY)^{\frac nd}
=\left(\sum_{d|n}d(1-q^{\frac nd})\textsf{v}_d^q(X)^{\frac nd}\right)\left(\sum_{d|n}d(1-q^{\frac nd})\textsf{v}_d^q(Y)^{\frac nd}\right)
\end{align*}
for all $n\ge 1.$
In a recursive way, one easily sees that
$$\textsf{v}_n^q(XY) \in \mathbb Q[q][\textsf{v}_d^q(X),\textsf{v}_d^q(Y): d|n].$$
Now, we are ready to prove that
$$\textsf{v}_n^q(XY) \in \mathbb Z[q][\textsf{v}_d^q(X),\textsf{v}_d^q(Y): d|n]$$ for all $n\ge 1.$
Let $p$ be any prime divisor of $n$.
Then
\begin{align*}
&n(1-q)\textsf{v}_n^q(XY)\\
&=(\sum_{d|n}d(1-q^{\frac nd})\textsf{v}_d^q(X)^{\frac nd})\left(\sum_{d|n}d(1-q^{\frac nd})\textsf{v}_d^q(Y)^{\frac nd}\right)
-\sum_{d|n \atop d<n}d(1-q^{\frac nd})\textsf{v}_d^q(XY)^{\frac nd}\\
&\equiv( \sum_{d| \frac np}d(1-(q^p)^{\frac {n}{pd}})(\textsf{v}_d^q(X)^p)^{\frac {n}{pd}})
\left(\sum_{d| \frac np}d(1-(q^p)^{\frac {n}{pd}})(\textsf{v}_d^q(Y)^p)^{\frac {n}{pd}}\right)
-\sum_{d|\frac{n}{p}}d(1-q^{\frac nd})(\textsf{v}_d^q(XY)^p)^{\frac {n}{pd}}\\
&\equiv \sum_{d| \frac np}d(1-(q^p)^{\frac {n}{pd}})
\left((\textsf{v}_d^q(XY)|_{{q \mapsto q^p \atop \textsf{v}_i^q(X)\mapsto (\textsf{v}_i^q(X))^p} \atop \textsf{v}_i^q(Y)\mapsto (\textsf{v}_i^q(Y))^p})^{\frac {n}{pd}}-(\textsf{v}_d^q(XY)^p)^{\frac {n}{pd}}\right)\\
&\equiv 0 \quad \text{(by Lemma \ref{key modular equivalence 1})(a)).}
\end{align*}
Here `$\equiv$' denotes `$\equiv \pmod {p^{\nu_p(n)}}$'.
Since $p$ is arbitrary, our assertion is verified.
\qed

\vskip 2mm
For a partition $\ld=(\ld_1, \ldots,\ld_l)$, let
\begin{align*}
&\textsf{u}_\ld^q(X):=\textsf{u}_{\ld_1}^q(X)\textsf{u}_{\ld_2}^q(X)\cdots \textsf{u}_{\ld_l}^q(X),\\
&\textsf{v}_\ld^q(X):=\textsf{v}_{\ld_1}^q(X)\textsf{v}_{\ld_2}^q(X)\cdots \textsf{v}_{\ld_l}^q(X),\\
&\textsf{h}^q_\ld(X):=\textsf{h}^q_{\ld_1}(X)\textsf{h}^q_{\ld_2}(X)\cdots \textsf{h}^q_{\ld_l}(X),\\
&\textsf{g}^q_\ld(X):=\textsf{g}^q_{\ld_1}(X)\textsf{g}^q_{\ld_2}(X)\cdots \textsf{g}^q_{\ld_l}(X).
\end{align*}
Let $X=\{x_i:i\ge 1 \}$ be an alphabet,
$\Lambda_{\mathbb Z}$ the ring of symmetric functions in infinitely many variables $x_i,(i\ge 1),$
over $\mathbb Z$, and $Par$ the set of all partitions.
For simplicity, let
$$\mathbb Z[q]_{1-q}:=\mathbb Z[q]\left[\frac {1}{1-q}\right].$$
We prove the following theorem.
\begin{thm} With the notation above, we have

{\rm (a)}$\{\textsf{u}_\ld^q(X): \ld \in Par\}$ is a $\mathbb Q[q]$-basis of $\mathbb Q[q]\otimes \Lambda_{\mathbb Z}$.

{\rm (b)}$\{\textsf{h}^q_\ld(X): \ld \in Par\}$ is a $\mathbb Q[q]$-basis of $\mathbb Q[q]\otimes \Lambda_{\mathbb Z}.$

{\rm (c)}$\{\textsf{v}_\ld^q(X): \ld \in Par\}$ is a $\mathbb Z[q]_{1-q}$-basis of $\mathbb Z[q]_{1-q}\otimes \Lambda_{\mathbb Z}.$

{\rm (d)}$\{\textsf{g}^q_\ld(X): \ld \in Par\}$ is a $\mathbb Z[q]_{1-q}$-basis of $\mathbb Z[q]_{1-q}\otimes \Lambda_{\mathbb Z}.$
\end{thm}
\pf

(a)
From the first equation of \eqref{rel between u and p}
it follows that
$$p_n(X)-n\textsf{u}_n^q(X) \in \mathbb Z[q][\textsf{u}_d^q(X): d|n, d<n]$$
for all $n\ge 1$.
Since $\{p_\ld(X): \ld \in Par \}$ is a $\mathbb Q$-basis of $\mathbb Q\otimes \Lambda_{\mathbb Z}$,
it follows that
$\{\textsf{u}_\ld^q(X): \ld \in Par\}$ is a spanning set of $\mathbb Q[q]\otimes \Lambda_{\mathbb Z}$.
To see that it is linearly independent,
let
$$a_{\ld^{(1)}}\textsf{u}^q_{\ld^{(1)}}+a_{\ld^{(2)}}\textsf{u}^q_{\ld^{(2)}}+\cdots+a_{\ld^{(r)}}\textsf{u}^q_{\ld^{(r)}}=0, \text{ with }a_{\ld^{(i)}}\in \mathbb Q[q],$$
where $\ld^{(1)}\unrhd \ld^{(2)} \unrhd \cdots \unrhd \ld^{(r)}$ in the dominance order.
Writing the left hand side in terms of $p_\ld(X), (\ld \in Par)$ yields that
$$a_{\ld^{(1)}}p_{\ld^{(1)}}+ \text{ lower terms} =0.$$
Hence, we can conclude that $a_{\ld^{(1)}}=0$.
Repeating this process, we can see that all the coefficients vanish.

(b)
Since
\begin{align*}
\textsf{h}^q_n(X)=\sum_{a_1+2a_2+\cdots+na_n=n \atop a_i \ge 0}
{\textsf{u}_1^q(X)}^{a_1}{\textsf{u}_2^q(X)}^{a_2}\cdots {\textsf{u}_n^q(X)}^{a_n},
\end{align*}
it follows that
$$\textsf{h}^q_n(X)-\textsf{u}_n^q(X) \in \mathbb Z[\textsf{u}_d^q(X): d|n, d<n]$$ for all positive integers $n$.
Thus our assertion follows from (a).

(c)
To begin with, we will show that
$$\textsf{v}_n^q(X)\in \mathbb Z[q]_{1-q}\otimes \Lambda_{\mathbb Z}, \,\, (n\ge 1).$$
From the identity
$$\prod_{n\ge 1}\frac{1-q\textsf{v}_n^q(X)t^n}{1-\textsf{v}_n^q(X)t^n}=H(X,t)$$
it follows that $h_n(X)$ can be expressed as
\begin{equation}\label{explict form of h interms of v}
\sum_{i+j=n \atop 0 \le i,j \le n}\left(\sum_{a_1+2a_2+\cdots+ia_i=i \atop a_l \ge 0,  (1\le l \le i)}
{\textsf{v}_1^q(X)}^{a_1}{\textsf{v}_2^q(X)}^{a_2}\cdots {\textsf{v}_n^q(X)}^{a_n}\right)
\left(\sum_{b_1+b_2+\cdots+b_s=j \atop b_k >0,  (1\le k \le s)}(-q)^{j}
{\textsf{v}_{b_1}^q(X)}\cdots {\textsf{v}_{b_s}^q(X)}\right).
\end{equation}
In case where $n=1$, our claim is obvious since
$$\textsf{v}_n^q(X)=\frac {1}{1-q}h_1(X).$$
When $n>1$, assume that our assertion is true for all positive integers less than $n$.
By \eqref{explict form of h interms of v} we can deduce that
$$h_n(X)-(1-q)\textsf{v}_n^q(X)\in \mathbb Z[q][\textsf{v}_i^q(X): 1\le i\le n-1],$$
and thus
$$\textsf{v}_n^q(X)-\frac {1}{1-q}h_n(X) \in \mathbb Z[q][\textsf{v}_i^q(X): 1\le i\le n-1].$$
Now, the claim is immediate from the induction hypothesis.
We are now ready to show that $\{\textsf{v}_\ld^q(X): \ld \in Par\}$ is a basis.
In view of \eqref{explict form of h interms of v},
$\{\textsf{v}_\ld^q(X): \ld \in Par\}$ is obviously a spanning set of $\mathbb Z[q]_{1-q}\otimes \Lambda_{\mathbb Z}.$
The linear independence can be shown in the same way as in (a).

(d) It can be verified in the same manner as in (b).
\qed

\begin{rem}{\rm
It should be noted that
$\{\textsf{u}_\ld^q(X): \ld \in Par\}$ is not a $\mathbb Z[q]$-basis of $\mathbb Z[q]\otimes \Lambda_{\mathbb Z}$.
To illustrate this, define symmetric functions $q_n(X), (n\ge 1)$ via the relations:
\begin{align*}
p_n(X)=\sum_{d|n}dq_d(X)^{\frac nd},\,\, (n\ge 1).
\end{align*}
For a partition $\ld=(\ld_1, \ldots,\ld_l)$, let
$$q_\ld(X):=q_{\ld_1}(X)q_{\ld_2}(X)\cdots q_{\ld_l}(X).$$
It is easy to see that
$\{q_\ld(X): \ld \in Par\}$ is a $\mathbb Z$-basis of $\Lambda_{\mathbb Z}$ (see \cite{R,ST}).
However, the system of equations
$$\sum_{d|n}d\left[\frac nd \right]_q\textsf{u}_d^q(X)^{\frac nd}=\sum_{d|n}dq_d(X)^{\frac nd},\,\, (n\ge 1),$$
imply that
$$q_p(X)=\textsf{u}_p^q(X)+\frac{[p]_q-1}{p}\textsf{u}_1^q(X)^{p},$$
where $p$ is a prime.
Obviously $\frac{[p]_q-1}{p}\notin \mathbb Z[q].$
}\end{rem}

\subsection {Connection between ${\overline {W}}^{q(q)}$ and ${\overline \Lambda}^{g(q)}$}
\label{connection W & Lambda}
As before, let $g(q)$ be an arbitrary polynomial in $\mathbb Z[q]$.
Given a $\mathbb Z[q]$-algebra $A$, consider the map $\Theta_A: {\overline {W}}^{g(q)}(A) \to {\overline \Lambda}^{g(q)}(A)$
defined by
$$(a_n)_{n\in \mathbb N} \mapsto \prod_{n\ge 1}\frac{1}{1-a_nt^n}.$$
Since the ring structure of ${\overline \Lambda}^{g(q)}(A)$ is transported from that of ${\overline {W}}^{g(q)}(A)$ via $\Theta_A$,
the assignment $A\mapsto \Theta_A$
induces a natural isomorphism $\Theta: {\overline {W}}^{g(q)} \to {\overline \Lambda}^{g(q)}$
such that the diagram
\begin{equation*}
\begin{CD}
{\overline {W}}^{g(q)} @>\Theta>\cong>{\overline \Lambda}^{g(q)} \\
@V{{\overline \Phi}^{g(q)}}VV  @V{\overline \Upsilon}^{g(q)}VV \\
\textsf{gh}\circ F^{\textsf{$\mathbb Z[q]$-Algebras}}_{\textsf{Rings}}
@>>{int}>\textsf{Gh}\circ F^{\textsf{$\mathbb Z[q]$-Algebras}}_{\textsf{Rings}}
\end{CD}
\end{equation*}
is commutative.
Here ${int}$ denotes the natural transformation such that
$int_A((a_n)_{n\in \mathbb N})=\sum_{n\ge 1}a_nt^{n-1}$
(refer to Eq.\eqref{How to define addition and multiplication g(q) case}).

In the rest of this subsection, we introduce some product identities which arise naturally in our context.
Let $x_n\, (n\ge 1)$ be indeterminates and assume that $g(q)$ ranges over the set $\mathbb Z[q]\setminus \{1\}$.
For $h(q)\in \mathbb Z[q]\setminus \{1\}$, consider the following system of equations:
\begin{equation}\label{explicit form of natural transformation g(q) case version}
\sum_{d|n}d(1-g(q)^{\frac nd})x_d^{\frac nd}=\sum_{d|n}d(1-h(q)^{\frac nd})\textsf{T}_d^{\frac nd}, \,\, (n\ge 1).
\end{equation}

\begin{lem} \label{To define natural transformation:polynpmial case}
Assume that $g(q)$ ranges over the set $\mathbb Z[q]\setminus \{1\}$.
Then, $\textsf{T}_n\in \mathbb Z[g(q)][x_d: d|n]$ for every positive integer $n$
if and only if $h(q)=0,2$.
\end{lem}
\pf
First, assume that $\textsf{T}_n\in \mathbb Z[g(q)][x_d: d|n]$ for every positive integer $n$.
From Eq.\eqref{explicit form of natural transformation g(q) case version} it follows that
$$\textsf{T}_n-\frac {1-g(q)}{1-h(q)}x_n \in \mathbb Z[g(q)][x_d: d|n, d<n]$$
for all $n\ge 1$.
Since $g(q)$ is arbitrary, it follows that $1-h(q)=\pm 1$, that is, $h(q)=0$ or $2$.

To prove the converse, assume that $h(q)=0$ or $2$.
We will prove the desired result by using mathematical induction on $n$.
When $n=1$, our assertion is trivial since $\textsf T_1=\pm (1-g(q))x_1$ by Eq.\eqref{explicit form of natural transformation g(q) case version}.
When $n>1$, assume that it holds for all positive integers less than $n$.
First, let $h(q)=0$.
Let $p$ be any prime divisor of $n$.
Then
\begin{align*}
n\textsf {T}_n
&=\sum_{d|n}d(1-g(q)^{\frac nd})x_d^{\frac nd}
-\sum_{d|n \atop d<n}d\textsf {T}_d^{\frac nd}\\
&\equiv \sum_{d| \frac np}d(1-(g(q)^p)^{\frac {n}{pd}})(x_d^p)^{\frac nd}
-\sum_{d| \frac np}d(\textsf {T}_d^p)^{\frac {n}{pd}} \pmod{p^{\nu_p(n)}}\\
&\equiv \sum_{d| \frac np}d((\textsf {T}_d|_{g(q) \mapsto g(q)^p \atop {x_i \mapsto x_i^p}})^{\frac {n}{pd}}
-(\textsf {T}_d^p)^{\frac {n}{pd}})\pmod{p^{\nu_p(n)}}\\
&\equiv 0 \pmod{p^{\nu_p(n)}} \quad \text{(by Lemma \ref{key modular equivalence 1}(a)).}
\end{align*}
On the other hand, if $h(q)=2$, then
\begin{align*}
-n\textsf {T}_n
&=\sum_{d|n}d(1-g(q)^{\frac nd})x_d^{\frac nd}
-\sum_{d|n \atop d<n}d(1-2^{\frac nd})\textsf {T}_d^{\frac nd}\\
&\equiv \sum_{d| \frac np}d(1-(g(q)^p)^{\frac {n}{pd}})(x_d^p)^{\frac nd}
-\sum_{d| \frac np}d(1-2^{\frac {n}{pd}})(\textsf {T}_d^p)^{\frac {n}{pd}} \pmod{p^{\nu_p(n)}}\\
&\equiv \sum_{d| \frac np}d(1-2^{\frac {n}{pd}})((\textsf {T}_d|_{g(q) \mapsto g(q)^p \atop {x_i \mapsto x_i^p}})^{\frac {n}{pd}}
-(\textsf {T}_d^p)^{\frac {n}{pd}}) \pmod{p^{\nu_p(n)}}\\
&\equiv 0 \pmod{p^{\nu_p(n)}} \quad \text{(by Lemma \ref{key modular equivalence 1}(a)).}
\end{align*}
Since $p$ is arbitrary, we can conclude that $\textsf T_n$ is divisible by $n$.
This completes the proof.
\qed

\vskip 2mm
Let $h(q)=0,2$.
For a $\mathbb Z[q]$-algebra $A$, let $\omega^{g(q),h(q)}_A:{\overline {W}}^{g(q)}(A) \to {\overline {W}}^{h(q)}(A)$ be the map defined by 
$$(a_n)_{n\in \mathbb N} \mapsto (\textsf{T}_n|_{x_d \mapsto a_d})_{n\in \mathbb N}.$$
Then the assignment $A\mapsto \omega^{g(q),h(q)}_A$ induces a natural transformation
$\omega^{g(q),h(q)}:{\overline {W}}^{g(q)} \to {\overline {W}}^{h(q)}$ such that the diagram
\begin{equation*}
\begin{CD}
{\overline {W}}^{g(q)} @>\omega^{g(q),h(q)}>>{\overline {W}}^{h(q)}@>\Theta>>{\overline \Lambda}^{h(q)}\\\
@V{{\overline \Phi}^{g(q)}}VV  @V{\overline \Phi}^{h(q)}VV @V{\overline \Upsilon}^{h(q)}VV\\
\textsf{gh}\circ F^{\textsf{$\mathbb Z[q]$-Algebras}}_{\textsf{Rings}}
@>>{\rm id}>\textsf{gh}\circ F^{\textsf{$\mathbb Z[q]$-Algebras}}_{\textsf{Rings}}
@>>{\it int}>\textsf{Gh}\circ F^{\textsf{$\mathbb Z[q]$-Algebras}}_{\textsf{Rings}}
\end{CD}
\end{equation*}
is commutative.

Suppose that $l,m$ are arbitrary integers.
Combining Eq.\eqref{explicit form of natural transformation g(q) case version} and Lemma \ref{To define natural transformation:polynpmial case}
yields the following proposition.

\begin{prop}
Let $A$ be any $\mathbb Z[q]$-algebra and assume that $g(q)$ ranges over $\mathbb Z[q]\setminus \{1\}$.
For $(a_n)_{n\in \mathbb N}\in {\overline {W}}^{g(q)}(A)$,
let
$(b_n)_{n\in \mathbb N}:=\omega^{g(q),0}_A((a_n)_{n\in \mathbb N})$ and
$(c_n)_{n\in \mathbb N}:=\omega^{g(q),2}_A((a_n)_{n\in \mathbb N}).$
Then we have
\begin{align*}
\prod_{n\ge 1}\frac{1-g(q)a_nt^n}{1-a_nt^n}
=\prod_{n\ge 1}\frac{1}{1-b_nt^n}
=\prod_{n\ge 1}\frac{1-2c_nt^n}{1-c_nt^n}.
\end{align*}
\end{prop}
\vskip 2mm
Define $\textsf{z}_n (n\ge 1)$ via the following system of equations:
\begin{equation*}\label{explicit form of natural transformation}
\sum_{d|n}d(1-m^{\frac nd})x_d^{\frac nd}=\sum_{d|n}d(1-l^{\frac nd})\textsf{z}_d^{\frac nd}, \quad  (n\ge 1).
\end{equation*}
Letting $g(q)=q$ and then specializing $q$ into an integer $m$ in the above paragraph,
we can derive the following corollary.

\begin{cor} Let $m$ vary the set of integers. Then we have the following.

{\rm (a)} With the above notation,
$\textsf{z}_n\in \mathbb Z[x_d: d|n]$ for every positive integer $n$
if and only if $l=0,2$

{\rm (b)}
Let $l=0,2$.
For a $\mathbb Z[q]$-algebra $A$, let
$$\bar \omega^{m,l}_A:W^{m}(A) \to W^{l}(A), \quad (a_n)_{n\in \mathbb N} \mapsto (\textsf{z}_n|_{x_d \mapsto a_d})_{n\in \mathbb N}.$$
Then the assignment $A\mapsto \omega^{m,l}_A$ induces a natural transformation
$\bar \omega^{m,l}:W^{m} \to W^{l}$ such that the diagram
\begin{equation*}
\begin{CD}
W^{\, m} @>\bar \omega^{m,l}>>W^{\, l}@>\Theta>>\Lambda^{l}\\\
@V{\Phi^{m}}VV  @V\Phi^{l}VV @V\Upsilon^{l}VV\\
\textsf{gh}
@>>{\rm id}>\textsf{gh}
@>>{\it int}>\textsf{Gh}
\end{CD}
\end{equation*}
is commutative.

{\rm (c)}
Let $A$ be any commutative ring with unity and let $m$ be any integer.
For $(a_n)_{n\in \mathbb N}\in W^{m}(A)$,
let
$(b_n)_{n\in \mathbb N}:=\bar \omega^{m,0}_A((a_n)_{n\in \mathbb N})$ and
$(c_n)_{n\in \mathbb N}:=\bar \omega^{m,2}_A((a_n)_{n\in \mathbb N}).$
Then we have
\begin{align*}
\prod_{n\ge 1}\frac{1-ma_nt^n}{1-a_nt^n}
=\prod_{n\ge 1}\frac{1}{1-b_nt^n}
=\prod_{n\ge 1}\frac{1-2c_nt^n}{1-c_nt^n}.
\end{align*}
\end{cor}

\section{$q$-deformation of necklace rings}
\label{q-Deformations of necklace rings and their specializations}
As explained in the introduction, classical necklace rings come as a partner of classical Witt rings.
Motivated by this phenomenon, in this section, we provide functors
$${{\overline {\textsf B}}}^{g(q)}:  {\textsf{$\mathbb Z[q]\otimes {\Psi}$-Rings}} \to \textsf{RINGS}$$
attached to each polynomial $g(q)\in \mathbb Z[q]$
and $${\textsf B}^m: \textsf{$\Psi$-Rings} \to \textsf{RINGS}$$ attached to each integer $m\in \mathbb Z.$
The relationship between ${{\overline {W}}}^{g(q)}$ and ${{\overline {\textsf B}}}^{g(q)}$ and that between ${W}^{m}$ and ${\textsf B}^{m}$
will be studied intensively in Section \ref{Relationship between Witt rings and necklace rings}.

\subsection{$q$-Deformation of necklace rings}
Throughout this subsection, we assume that $\mathbb Z[q]$ is equipped with the $\Psi$-ring structure such that
$\Psi^n(f(q))=f(q^n)$
for all $f(q) \in \mathbb Z[q]$ and $n\ge 1$.
Let $\textsf E$ be the $\Psi$-ring freely generated by $x_n,y_n\,(n\ge 1)$.
It is easy to see that ${\textsf F}:=\mathbb Z[q]\otimes {\textsf E}$ is made into a $\Psi$-ring if we define
$$\Psi^n(a\otimes b)=\Psi^n(a)\otimes \Psi^n(b),\, a\in \mathbb Z[q],\,  b\in {\textsf E}.$$

\begin{rem}{\rm
Proposition \ref {psi=ld} implies that
the polynomial ring $\mathbb Z[q]$ can be equipped with the $\ld$-ring structure
with the Adams operations defined as above.
}\end{rem}
\vskip 2mm

Let ${\textsf R}:=\mathbb Q[q]\otimes {\textsf E}$.
Given a polynomial $g(q)\in \mathbb Z[q]$,
consider the map
$${\overline \varphi}^{g(q)}_{\textsf R}: {\textsf R}^{\mathbb N} \to {\textsf R}^{\mathbb N}, \quad
(a_n)_{n\in \mathbb N} \mapsto \left(\sum_{d|n}d(1-g(q)^{\frac nd}) \Psi^{\frac nd}(a_d)\right)_{n\in \mathbb N}.$$
Let $X=(x_n)_{n\in \mathbb N}$ and $Y=(y_n)_{n\in \mathbb N}$.
Let $P=(P_n)_{n\in \mathbb N}$
be the element determined by the relation:
\begin{equation}\label{multiplication defining poly}
\sum_{d|n}d\left[\frac nd\right]_{g(q)}\Psi^{\frac nd}(P_d)
=(1-g(q))\left(\sum_{d|n}d\left[\frac nd\right]_{g(q)}\Psi^{\frac nd}(x_d)\right)
\left(\sum_{d|n}d \left[\frac nd\right]_{g(q)}\Psi^{\frac nd}(y_d)\right),\, (n\ge 1).
\end{equation}
Multiplying either side of Eq.\eqref{multiplication defining poly} by $1-g(q)$ yields that
$${\overline \varphi}^{g(q)}_{{\textsf R}}(P)=
{\overline \varphi}^{g(q)}_{\textsf R}(X)\,\,{\overline \varphi}^{g(q)}_{\textsf R}(Y).$$

\begin{lem}\label{key lemma in the proof of existence of q deform of necklace ring}
Let $n$ be an arbitrary positive integer.
Then
$P_n \in \mathbb Z[q][\Psi^{\frac nd}(x_d),\Psi^{\frac nd}(y_d):d|n].$
\end{lem}

\pf
If $g(q)=1$, then $P_n=0$ for all $n\ge 1$.
Indeed it can be verified by applying the mathematical induction on $n$ to Eq.\eqref{multiplication defining poly}.
So, from now on, assume $g(q)\ne 1$.
Note $P_1=(1-g(q))x_1y_1$.
If $n>1$, assume that the desired assertion holds for all positive integers less than $n$.
Multiplying either side of Eq.\eqref{multiplication defining poly} by $1-g(q)$, we have
\begin{equation}\label{necklace functor corr to g(q) multiplied }
\begin{aligned}
&n(1-g(q))P_n\\
&=\left(\sum_{d|n}d(1-g(q)^{\frac nd})\Psi^{\frac nd}(x_d)\right)
\left(\sum_{d|n}d (1-g(q)^{\frac nd}) \Psi^{\frac nd}(y_d)\right)
-\sum_{d|n \atop d<n}d(1-g(q)^{\frac nd})\Psi^{\frac nd}(P_d).
\end{aligned}
\end{equation}
Due to the induction hypothesis, it is clear that the right hand side of the above equation is contained in $\mathbb Z[q][\Psi^{\frac nd}(x_d),\Psi^{\frac nd}(y_d):d|n]$.
The proof of Lemma \ref{key lemma in the proof of existence of q deformation of Witt rings}
shows that our claim is almost straightforward if one can
prove that the right hand side of Eq.\eqref{classification of Witt rings; defining polynomial} is divisible by
$n$.
Let $p$ be a prime divisor of $n$.
Then we have the following congruence:
\begin{equation}\label{necklace functor corr to g(q) reduced}
\begin{aligned}
n(1-g(q))P_n
&\equiv  \left(\sum_{d|\frac np}d(1-g(q)^{\frac nd})\Psi^{\frac nd}(x_d)\right)
\left(\sum_{d|\frac np}d (1-g(q)^{\frac nd}) \Psi^{\frac nd}(y_d)\right)\\
&\quad -\sum_{d|\frac np}d(1-g(q)^{\frac nd})\Psi^{\frac nd}(P_d)\pmod{p^{\nu_p(n)}}.
\end{aligned}
\end{equation}
In view of Lemma \ref{key modular equivalence 1}(a),
one has
\begin{equation}\label{key modular relation}
\Psi^p(1-g(q)^{\frac {n}{pd}})=1-g(q^p)^{\frac {n}{pd}}\equiv 1-g(q)^{\frac {n}{d}} \pmod {p^{\nu_p(\frac nd)}}
\end{equation}
for every divisor $d$ of $\frac np$.
Therefore,
\begin{equation*}
\begin{aligned}
&\left(\sum_{d|\frac np}d(1-g(q)^{\frac nd})\Psi^{\frac nd}(x_d)\right)
\left(\sum_{d|\frac np}d (1-g(q)^{\frac nd}) \Psi^{\frac nd}(y_d)\right)\\
&\equiv \Psi^p \left(\left(\sum_{d|\frac np}d(1-g(q)^{\frac {n}{pd}})\Psi^{\frac {n}{pd}}(x_d)\right)
\left(\sum_{d|\frac np}d (1-g(q)^{\frac {n}{pd}}) \Psi^{\frac {n}{pd}}(y_d)\right)\right)\pmod{p^{\nu_p(n)}}\\
&=\Psi^p \left(\sum_{d|\frac np}d(1-g(q)^{\frac {n}{pd}})\Psi^{\frac {n}{pd}}(P_d)\right)
\quad \text{(by Eq.\eqref{multiplication defining poly})}\\
&= \sum_{d|\frac np}d(1-g(q)^{\frac {n}{d}})\Psi^{\frac {n}{d}}(P_d).
\end{aligned}
\end{equation*}
Apply this to \eqref{necklace functor corr to g(q) reduced} to get
$$n(1-g(q))P_n\equiv 0 \pmod{p^{\nu_p(n)}}.$$
Consequently $n(1-g(q))P_n\equiv 0 \pmod{n}$, as required.
\qed

\vskip 2mm
Let us define the category {\textsf{$\mathbb Z[q]\otimes {\Psi}$-Rings}}.
Its objects consist of the $\mathbb Z[q]$-algebras that are equipped with
the structure of a $\Psi$-ring such that $\Psi^n(q)=q^n$ for all $n\ge 1.$
The morphisms are $\mathbb Z[q]$-algebra homomorphisms compatible with $\Psi$-operations.
By virtue of Lemma \ref{key lemma in the proof of existence of q deform of necklace ring}, one can construct a functor
$${{\overline {\textsf B}}}^{g(q)}: {\textsf{$\mathbb Z[q]\otimes {\Psi}$-Rings}}\to {\textsf{RINGS}}$$
such that, for each object $R$,
the addition ${{\overline {\textsf B}}}^{g(q)}(A)$ is defined componentwise and
the multiplication is defined via the universal polynomials $P_n$'s $(n\ge 1)$
(refer to the proof of Theorem \ref{Construction of new Witt rings attached to polynomials}).
The following theorem can be derived essentially in the same way as in Theorem \ref{Construction of new Witt rings attached to polynomials}.
\begin{thm}\label{defining B(g(q))}
Given an polynomial $g(g)\in \mathbb Z[q]$, there exists a functor
$${{\overline {\textsf B}}}^{g(q)}: {\textsf{$\mathbb Z[q]\otimes {\Psi}$-Rings}}\to {\textsf{RINGS}}$$
subject to the following conditions$:$
\begin{enumerate}
\item
As a set, ${{\overline {\textsf B}}}^{g(g)}(R)=R^\mathbb N$.
\item
For every morphism $f:R\to S$ and every $\alpha \in {{\overline {\textsf B}}}^{g(g)}(R)$ one has
${{\overline {\textsf B}}}^{g(g)}(f)(\alpha)=f \circ \alpha$.
\item
The map ${\overline \varphi}^{g(g)}_{R}:{{\overline {\textsf B}}}^{g(g)}(R)\to \textsf {gh}(R)$ defined by
$$(a_n)_{n\in \mathbb N} \mapsto \left(\sum_{d|n}d(1-g(q)^{\frac nd})\Psi^{\frac nd}(a_d)\right)_{n\in \mathbb N}$$
is a ring homomorphism.
\end{enumerate}
If $g(q)\ne 1$, then ${{\overline {\textsf B}}}^{g(g)}$ is uniquely determined by the conditions (1)-(3).
\end{thm}

\begin{rem}{\rm
(a) In the above theorem, one can see that ${{\overline \varphi}}^{g(q)}$ is a natural transformation.

(b) Obviously $\mathbb Z[q]$ is equipped with the $\Psi$-ring structure with
$\Psi^n={\rm id}$ for all $n\ge 1$. With this $\Psi$-ring structure, define
$P_n$ as in Eq.\eqref{multiplication defining poly}.
But, in this case,
$$\textsf {P}_n \notin \mathbb Z[q][\Psi^{\frac nd}(x_d),\Psi^{\frac nd}(y_d):d|n]$$
in general.
For instance, if $g(q)=q$, then
$$\textsf {P}_2=(1-q^2)\Psi^2(x_1)y_2+(1-q^2)x_2\Psi^2(y_1)+2(1-q)x_2y_2+\frac{q-q^3}{2}\Psi^2(x_1y_1).$$
}\end{rem}

\subsection{Necklace ring attached to an integer}
In this subsection, let ${\textsf F}$ will denote the $\Psi$-ring freely generated by $x_n,y_n\,(n\ge 1)$.
As a ring, it is the polynomial ring
$$\mathbb Z[\Psi^k(x_n),\Psi^k(y_n): k,n\ge 1].$$
Let ${\textsf R}:=\mathbb Q[q]\otimes {\textsf F}$.
Given an integer $m$, let us consider the map
$$\varphi^m_{\textsf R}: \textsf R^{\mathbb N} \to \textsf R^{\mathbb N},\quad
(a_n)_{n\in \mathbb N}\mapsto
\left(\sum_{d|n}d (1-m^{\frac {n}{d}})\Psi^{\frac nd}(a_d)\right)_{n\in \mathbb N}.$$
Let $X=(x_n)_{n\in \mathbb N}$ and $Y=(y_n)_{n\in \mathbb N}$, and define
$\textsf P=(\textsf P_n)_{n\in \mathbb N}$ recursively via the following relations:
\begin{equation}\label{new q deform of necklace ring}
\sum_{d|n}d \left[\frac nd \right]_m \Psi^{\frac nd}(\textsf {P}_d)
=(1-m)\left(\sum_{d|n}d \left[\frac nd \right]_m \Psi^{\frac nd}(x_d)\right)
\left(\sum_{d|n}d \left[\frac nd \right]_m\Psi^{\frac nd}(y_d)\right),\quad (n\ge 1),
\end{equation}
equivalently,
$$\varphi^m_{\textsf R}(\textsf {P})=\varphi^m_{\textsf R}(X)\varphi^m_{\textsf R}(Y).$$
Exploiting the mathematical induction on $n$ one easily sees that
$$\textsf {P}_n \in \mathbb Q[\Psi^{\frac nd}(x_d),\Psi^{\frac nd}(y_d):d|n],  \quad (n\ge 1).$$
The proof of the following theorem can be done by a a slight modification of that of Lemma \ref{key lemma in the proof of existence of q deform of necklace ring}.

\begin{lem}\label{key lemma in the proof of existence of q deform of necklace ring 1}
Let $n$ be an arbitrary positive integer.
Then
$\textsf {P}_n \in \mathbb Z[\Psi^{\frac nd}(x_d),\Psi^{\frac nd}(y_d):d|n].$
\end{lem}

\vskip 2mm
By virtue of Lemma \ref{key lemma in the proof of existence of q deform of necklace ring 1},
one can construct a functor
$${{\textsf B}}^m:\textsf{$\Psi$-Rings} \to \textsf{RINGS}$$
such that the addition is defined componentwise and
the multiplication is defined via the universal polynomials $\textsf {P}_n$'s $(n\ge 1)$.
We prove the following theorem in the same way as in Theorem \ref{Construction of new Witt rings attached to polynomials}.
\begin{thm}\label{constructing B(m)}
Given an integer $m$, the functor ${{\textsf B}}^m$
satisfies the following conditions$:$
\begin{enumerate}
\item
As a set, ${{\textsf B}}^m(R)=R^\mathbb N$.
\item
For every $\Psi$-ring homomorphism $f:R\to S$
and every $\alpha \in {{\textsf B}}^m(R)$ one has
${{\textsf B}}^m(f)(\alpha)=f \circ \alpha$.
\item
The map
$\varphi^m_{R}:{{\textsf B}}^m(R)\to \textsf{gh}(R)$ defined by
$$(a_n)_{n\in \mathbb N} \mapsto
\left(\sum_{d|n}d(1-m^{\frac nd})\Psi^{\frac nd}(a_d)\right)_{n\in \mathbb N}$$
is a ring homomorphism.
\end{enumerate}
If $m\ne 1$, then it is uniquely determined by (1)-(3).
\end{thm}

\begin{rem}{\rm
(a) Theorem \ref{constructing B(m)}
cannot be deduced from Theorem \ref{defining B(g(q))} by the specialization $q \mapsto m$
since it is not compatible with $\Psi^n,\, (n\ge 1).$
For instance,
$$(1-q^n)\Psi^n(x) \stackrel {q\mapsto m}{\longrightarrow} (1-m^n)\Psi^n(x),$$
whereas
$$\Psi^n((1-q)x)  \stackrel {q\mapsto m}{\longrightarrow} \Psi^n((1-m)x)=(1-m)\Psi^n(x).$$

(b) The functors ${\overline {\textsf B}}^{g(q)}$ and $\textsf B^m$
have the exactly same classification theorems as in Section \ref{Classification Theorem( Witt ring)},
all of which can be obtained by the slight modification of the arguments in Section \ref{Classification Theorem( Witt ring)}.
}\end{rem}

\subsection{Structure constants of ${\overline {\textsf B}}^{g(q)}$ and $\textsf B^m$}
\label{Structure constants}

\subsubsection{}
In the previous section, we have shown that
$${\overline {\textsf B}}^{g(q)}(X){\overline {\textsf B}}^{g(q)}(Y)=(P_n),$$
where
$$P_n=\sum_{i,j|n}d(i,j)\Psi^{\frac ni}(x_i)\Psi^{\frac nj}(y_j).$$
Here we focus on the explicit form of $d(i,j)$.
To begin with, let us introduce definitions and notation necessary to develop our arguments.
For a positive integer $n$, let $D_n$ be the poset consisting of all divisors of $n$
where $i \preceq j$ in $D_n$ if $j$ is divisible by $i$.
For positive integers $a,b\in D_n$ with $a \preceq b$, let
$\mathcal C(a,b)$ the set of all chains from $a$ to $b$ in $D_n$.
It is easy to see that it does not depend on the choice of $n$.

For positive integers $a,b\in D_n$ with $a \preceq b$ and
for any chain
$C=(a=x_r \prec \cdots \prec x_0=b)\in \mathcal C(a,b),$
let
$${\rm Wt}_{g(q)}(C):=(-1)^r\left[\frac{x_0}{x_1}\right]_{g(q)}\left[\frac{x_1}{x_2}\right]_{g(q^{\frac{x_0}{x_1}})}\cdots \left[\frac{x_{r-1}}{x_r}\right]_{g(q^{\frac{x_0}{x_{r-1}}})}.$$
and
$$\hat \mu_{g(q)}(b,a):=\sum_{C\in \mathcal C(a,b)}{\rm Wt}_{g(q)}(C).$$
It can be easily seen that $\hat \mu_{g(q)}(b,a)=\hat \mu_{g(q)}(b',a')$ if $b/a={b'}/{a'}$.
For each positive integer $n$, let us define $\hat \mu_{g(q)}(n)$ by $\hat \mu_{g(q)}(n,1)$.

\begin{prop} \label{explicit form of coefficients g(q) case}
Let $n$ be a positive integer and $i,j$ divisors of $n$.
Then we have
$$d(i,j)=\dfrac {ij}{n}\sum_{e|\frac{n}{[i,j]}}(1-g(q^{e}))\left[\frac{n}{ei}\right]_{g(q^e)}\left[\frac{n}{ej}\right]_{g(q^e)}\hat \mu_{g(q)}(e).$$
\end{prop}
\pf
Note that $n\textsf {P}_n$ equals
\begin{align}\label{explict form of p(n) 1  g(q) case}
(1-g(q))\left(\sum_{d_1\preceq n }d_1 \left[\frac {n}{d_1} \right]_{g(q)} \Psi^{\frac {n}{d_1}}(x_{d_1})\right)
\left(\sum_{d'_1\preceq n }d'_1 \left[\frac {n}{d'_1} \right]_{g(q)} \Psi^{\frac {n}{d'_1}}(y_{d'_1})\right)
-\sum_{d_1\prec n} \left[\frac {n}{d_1} \right]_{g(q)} \Psi^{\frac {n}{d_1}}(d_1 \textsf {P}_{d_1}).
\end{align}
If $d_1 \textsf {P}_{d_1}$ is replaced by
\begin{equation}\label{explict form of p(n) 2 g(q) case}
(1-{g(q)})\left(\sum_{d_2\preceq d_1 }d_2 \left[\frac {d_1}{d_2} \right]_{g(q)} \Psi^{\frac {d_1}{d_2}}(x_{d_2})\right)
\left(\sum_{d'_2\preceq d_1 }d'_2 \left[\frac {d_1}{d'_2} \right]_{g(q)} \Psi^{\frac {d_1}{d'_2}}(y_{d'_2})\right)
-\sum_{d_2\prec d_1} \left[\frac {d_1}{d_2} \right]_{g(q)} \Psi^{\frac {d_1}{d_2}}(d_2 \textsf {P}_{d_2}),
\end{equation}
then $n\textsf {P}_n$ is expressed as
\begin{align*}
&(1-g(q))\left(\sum_{d_1,d'_1\preceq n }d_1d'_1 \left[\frac {n}{d_1} \right]_{g(q)} \left[\frac {n}{d'_1} \right]_{g(q)}
\Psi^{\frac {n}{d_1}}(x_{d_1})
\Psi^{\frac {n}{d'_1}}(y_{d'_1})\right)\\
&-(1-{g(q^{\frac {n}{d_1}})})\left(\sum_{d_2,d'_2 \preceq d_1 \prec n }d_2d'_2 \left[\frac {n}{d_1} \right]_{g(q)}
\left[\frac {d_1}{d_2} \right]_{g(q^{\frac {n}{d_1}})}
\left[\frac {d_1}{d'_2} \right]_{g(q^{\frac {n}{d_1}})}
\Psi^{\frac {d_1}{d_2}}(x_{d_2})
\Psi^{\frac {d_1}{d'_2}}(y_{d'_2})\right)\\
&-\sum_{d_2\prec d_1 \prec n}  \left[\frac {n}{d_1} \right]_{g(q)}\left[\frac {d_1}{d_2} \right]_{g(q^{\frac {n}{d_1}})}
\Psi^{\frac {n}{d_2}}(d_2 \textsf {P}_{d_2}).
\end{align*}
Continue this process.
More precisely, replacing $d_i \textsf{P}_{d_i}\,(i\ge 2)$ by
\begin{equation}\label{explict form of p(n) 3 g(q) case}
\begin{aligned}
&(1-{g(q)})\left(\sum_{d_{i+1}\preceq d_i }d_{i+1} \left[\frac {d_i}{d_{i+1}} \right]_{g(q)} \Psi^{\frac {d_i}{d_{i+1}}}(x_{d_{i+1}})\right)
\left(\sum_{d'_{i+1}\preceq d_i}d'_{i+1} \left[\frac {d_i}{d'_{i+1}} \right]_{g(q)} \Psi^{\frac {d_i}{d'_{i+1}}}(y_{d'_{i+1}})\right)\\
&-\sum_{d_{i+1}\prec d_i } \left[\frac {d_i}{d_{i+1}} \right]_{g(q)} \Psi^{\frac {d_i}{d_{i+1}}}(d_{i+1} \textsf {P}_{d_{i+1}})
\end{aligned}
\end{equation}
in each step,
we finally arrive at the identity that
$$n\textsf {P}_n=\sum_{i,j|n}d(i,j)\Psi^{\frac ni}(x_i)\Psi^{\frac nj}(y_j),$$
where $d(i,j)$ is equal to
$$ij\sum_{r\ge 0 \atop i,j \preceq  d_{r}\prec d_{r-1}\prec \cdots \prec d_1\prec d_0:=n}(-1)^{r}(1-{g(q^{\frac {n}{d_{r}}})})
\left(\prod_{l=0}^{r-1}\left[\frac{d_{l}}{d_{l+1}}\right]_{g(q^{\frac {n}{d_l}})}\right) \left[\frac{d_{r}}{i}\right]_{g(q^{\frac {n}{d_r}})}\left[\frac{d_{r}}{j}\right]_{g(q^{\frac {n}{d_r}})},$$
where  $d_0:=n$ and $\prod_{l=0}^{r-1}\left[\frac{d_{l}}{d_{l+1}}\right]_{g(q^{\frac {n}{d_l}})}$ is set to be 1 when $r=0$.
Observe that the sum is over the set $\mathcal A(n;i,j)$ of all subposets of $D_n$ of the form
$\{d_0, d_1, \ldots, d_{r-1}, d_{r}\}$, where $r\ge 0$ and
$$i,j\preceq d_{r}\prec d_{r-1}\prec \cdots \prec d_2 \prec  d_1 \prec n.$$
Since $d_r$ can take any divisor with $[i,j]|d_r |n$, we have
$$\mathcal A(n;i,j)=\bigsqcup_{[i,j]|k|n}\mathcal C(k,n),$$
It is easy to see that the map
$$\eta:\mathcal A(n;i,j) \to \bigsqcup_{e|\frac {n}{[i,j]}}\mathcal C(e,\frac{n}{[i,j]}),\quad
\{d_0, d_1, \ldots, d_{r-1}, d_{r}\} \mapsto \left(\frac{d_{r-1}}{d_{r}},\ldots,\frac{d_0}{d_r}\right).$$
is bijective and weight-preserving, that is,
$${\rm Wt}_{g(q)}(C)={\rm Wt}_{g(q)}(\eta(C)), \quad C\in \mathcal A(n;i,j).$$
Consequently, we have
\begin{align*}
&\sum_{r\ge 0 \atop i,j \preceq  d_{r}\prec d_{r-1}\prec \cdots \prec d_1\prec d_0:=n}(-1)^{r}(1-{g(q^{\frac {n}{d_{r}}})})
\left(\prod_{l=0}^{r-1}\left[\frac{d_{l}}{d_{l+1}}\right]_{g(q^{\frac {n}{d_l}})}\right)
\left[\frac{d_{r}}{i}\right]_{g(q^{\frac {n}{d_r}})}\left[\frac{d_{r}}{j}\right]_{g(q^{\frac {n}{d_r}})}\\
&=\sum_{[i,j]|k|n}(1-g(q^{\frac nk}))\left(\sum_{r\ge 0 \atop {k=d_r \prec d_{r-1}\prec \cdots \prec d_2 \prec d_1 \prec d_0}}
(-1)^{r}\left(\prod_{l=0}^{r-1}\left[\frac{d_{l}}{d_{l+1}}\right]_{g(q^{\frac {n}{d_l}})}\right) \right)
\left[\frac{k}{i}\right]_{g(q^{\frac {n}{k}})}\left[\frac{k}{j}\right]_{g(q^{\frac {n}{k}})}\\
&=\sum_{[i,j]|k|n}(1-g(q^{\frac nk}))
\hat \mu_{g(q)}(n,k)
\left[\frac{k}{i}\right]_{g(q^{\frac {n}{k}})}\left[\frac{k}{j}\right]_{g(q^{\frac {n}{k}})}\\
&=\sum_{e|\frac {n}{[i,j]}}(1-g(q^e))\left[\frac {n}{ei}\right]_{g(q^e)}
\left[\frac {n}{ej}\right]_{g(q^e)}\hat \mu_{g(q)}(e)
\quad \text{ (by replacing $n/k$ by $e$)}.
\end{align*}
This completes the proof.
\qed

\begin{exm}{\rm
In case where $g(q)=q$, it can be easily seen that
$\hat \mu_{q}(e)=\left[e\right]_{q}\mu (e )$.
Therefore, we have
$$d(i,j)=\begin{cases}\dfrac {ij(1-q)}{n}\left[\frac{n}{i}\right]_{q}\left[\frac{n}{j}\right]_{q} & \text{ if }[i,j]=n,\\
0 &\text{ otherwise.}
\end{cases}$$
}\end{exm}

\vskip 2mm
Let
$\textsf B^m(X)\textsf B^m(Y)=(\textsf {P}_n)$, where
$\textsf {P}_n$ is of the form
$$\sum_{i,j|n}c(i,j)\Psi^{\frac ni}(x_i)\Psi^{\frac nj}(y_j).$$
For any chain
$C=(a=x_r \prec \cdots \prec x_0=b)\in \mathcal C(a,b),$
let $${\rm wt}_m(C):=
\begin{cases}(-1)^r\left[\frac{x_0}{x_1}\right]_m\cdots \left[\frac{x_{r-1}}{x_r}\right]_m& \text{ if }a<b\\
1&\text{ if } a=b\end{cases} $$
and
$$\mu_m(b,a):=\sum_{C\in \mathcal C(a,b)}{\rm wt}_m(C).$$
It can be easily seen that $\mu_m(b,a)=\mu_m(b',a')$ if $b/a={b'}/{a'}$.
For each positive integer $n$, let us define $\mu_m(n)$ by $\mu_m(n,1)$.
By the slight modification of the proof of Proposition \ref{explicit form of coefficients g(q) case}
we can derive the explicit form of $c(i,j)$.
\begin{cor} \label{explicit form of coefficients}
Let $n$ be a positive integer and $i,j$ divisors of $n$. Then we have
$$c(i,j)=\frac {(1-m)ij}{n}\sum_{e|\frac{n}{[i,j]}}\left[\frac{n}{ei}\right]_m \left[\frac{n}{ej}\right]_m\mu_m(e).$$
\end{cor}

\begin{exm}\label{concrete form of coefficients}{\rm

By simple calculation, one can show that
\begin{align*}
\textsf P_n
=\begin{cases}
\sum_{[i,j]=n}(i,j)\Psi^{\frac ni}(x_i)\Psi^{\frac nj}(y_j)& \text{ when } m=0,\\
0& \text{ when } m=1,\\
2\sum_{[i,j]=n \atop \frac ni, \frac nj \text{ are odd}}(i,j)\Psi^{\frac ni}(x_i)\Psi^{\frac nj}(y_j)&\text{ when } m=-1.
\end{cases}
\end{align*}
}\end{exm}

\begin{rem}{\rm
Corollary \ref{explicit form of coefficients} implies that
$$(1-m)ij\sum_{e|\frac{n}{[i,j]}}\left[\frac{n}{ei}\right]_m \left[\frac{n}{ej}\right]_m\mu_m(e)$$
is divisible by $n$.
It would be very nice to give a direct proof for this divisibility.
For instance, if $i=j=1$, then
$$(1-m)\sum_{e|n}{\left[\frac ne \right]^2_m} \mu_m(e)$$
is divisible by $n$.
}\end{rem}

\section{Functorial and structural properties}
\label{Functorial and structural properties}

\subsection{Induction}
Assume that ${\mathcal F}$ denotes
$\textsf{gh}\circ F^{\textsf{{$\mathbb Z[q]$-Algebras}}}_{\textsf{Rings}},$
$\textsf{gh}\circ F^{{\textsf{$\mathbb Z[q]\otimes {\Psi}$-Rings}}}_{\textsf{Rings}}$,
or $\textsf{gh}.$
For each $r\ge 1$ and each object $A$, consider the map
$${\textsf {Ind}}_{r,A}:\mathcal F(A)\to \mathcal F(A), \quad (a_n)_{n\in \mathbb N}\mapsto (ra_{\frac nr})_{n\in \mathbb N},$$
where
$a_{\frac nr}$ is set to be zero if $\frac nr \notin \mathbb N.$
It should be noted that this map is additive, but not multiplicative.
It is easy to show that the assignment $A\mapsto {\textsf {Ind}}_{r,A}$ induces a natural transformation
${\textsf {Ind}}_{r}:F^{\textsf{RINGS}}_{\textsf {Abel}}\circ \mathcal F \to F^{\textsf{RINGS}}_{\textsf {Abel}}\circ \mathcal F.$

Now assume that ${\mathcal G}$ denotes one of the following functors:
$${\overline {W}}^{g(q)},\,{\overline {\textsf B}}^{g(q)},\,
W^m, \text{ and }{{\textsf B}}^m, \quad \text{ where } g(q)\in \mathbb Z[q] \text{ and } m\ge 1.$$
For each $r\ge 1$ and each object $A$, consider the map
$${\textsl {Ind}}_{r,A}:\mathcal G(A)\to \mathcal G(A), \,(a_n)_{n\in \mathbb N}\mapsto (a_{\frac nr} )_{n\in \mathbb N}.$$

\begin{prop}\label{uniquness of induction B}
Suppose that $g(q)\in \mathbb Z[q]$ and $m\in \mathbb Z.$
For an arbitrary positive integer $r$ we have the following.

{\rm (a)} The assignment $A\mapsto {\textsl {Ind}}_{r,A}$ induces a natural transformation
$${\textsl {Ind}}_{r}:F^{\textsf{RINGS}}_{\textsf {Abel}}\circ {\overline {W}}^{g(q)} \to F^{\textsf{RINGS}}_{\textsf {Abel}}\circ {\overline {W}}^{g(q)}$$
such that
${\overline \Phi}^{g(q)}_A\circ {\textsl {Ind}}_{r,A}={\textsl {Ind}}_{r,A}\circ {\overline \Phi}^{g(q)}_A$ for every $\mathbb Z[q]$-algebra $A$.
If $g(q)\ne 1$, then it is uniquely determined by this commutativity.

{\rm (b)} The assignment $A\mapsto {\textsl {Ind}}_{r,A}$ induces a natural transformation
$${\textsl {Ind}}_{r}:F^{\textsf{RINGS}}_{\textsf {Abel}}\circ {\overline {\textsf B}}^{g(q)} \to F^{\textsf{RINGS}}_{\textsf {Abel}}\circ{\overline {\textsf B}}^{g(q)}$$
such that
${\overline \varphi}^{g(q)}_A\circ {\textsl {ind}}_{r,A}={\textsl {ind}}_{r,A}\circ {\overline \varphi}^{g(q)}_A$
for every object of {\textsf{$\mathbb Z[q]\otimes {\Psi}$-Rings}}.
If $m\ne 1$, then it is uniquely determined by this commutativity.

{\rm (c)} The assignment $A\mapsto {\textsl {Ind}}_{r,A}$ induces a natural transformation
$${\textsl {Ind}}_{r}:F^{\textsf{RINGS}}_{\textsf {Abel}}\circ W^m \to W^m$$
such that
$\Phi^m_A\circ {\textsl {Ind}}_{r,A}={\textsl {Ind}}_{r,A}\circ \Phi^m_A$ for every commutative ring $A$ with unity.
If $m\ne 1$, then it is uniquely determined by this commutativity.

{\rm (d)} The assignment $A\mapsto {\textsl {Ind}}_{r,A}$ induces a natural transformation
$${\textsl {Ind}}_{r}:F^{\textsf{RINGS}}_{\textsf {Abel}}\circ\textsf B^m \to \textsf B^m$$
satisfying
$\varphi^m_A\circ {\textsl {Ind}}_{r,A}={\textsl {Ind}}_{r,A}\circ \varphi^m_A$
for every $\Psi$-ring $A$.
If $m\ne 1$, then it is uniquely determined by this commutativity.
\end{prop}

\pf
Since all the assertions can be proved in the same way, we here prove only (a).
The proof will be accomplished in the following steps:

{\bf Step 1:} First, we will show that $\Phi_A^{g(q)}\circ {\textsl {Ind}}_{r,A}={\textsl {Ind}}_{r,A}\circ \Phi_A^{g(q)}$
for every $\mathbb Z[q]$-algebra $A$.
To do this it suffices to see that $(\Phi_{\rm ab})^{g(q)}_{A}\circ {\textsl {Ind}}_{r,A}={\textsl {Ind}}_{r,A}\circ (\Phi_{\rm ab})^{g(q)}_{A}$.
This is obvious since
\begin{align*}
(\Phi_{\rm ab})^{g(q)}_{A} \circ {\textsl{Ind}}_{r,A}((a_n)_{n\in \mathbb N})
&=(\Phi_{\rm ab})^{g(q)}_{A} ((a_{\frac nr})_{n\in \mathbb N})\\
&=\left(\sum_{d|n}d \left[\frac nd \right]_{g(q)}a_{\frac dr}^{\frac nd}\right)_{n\in \mathbb N}\\
&=\left(r\sum_{e|\frac nr}e\left[\frac {n}{er} \right]_{g(q)}a_{e}^{\frac {n}{er}}\right)_{n\in \mathbb N}
\quad \text{(by letting $e=d/r$)}\\
&={\textsf {Ind}}_{r,A}\circ (\Phi_{\rm ab})^{g(q)}_{A}((a_n)_{n\in \mathbb N}).
\end{align*}

{\bf Step 2:}
Second, we show that ${\textsl {Ind}}_{r,A}$ is additive for any $\mathbb Z[q]$-algebra $A$.
Let $\textsf F=\mathbb Z[q][x_n,y_n: n\ge 1]$, $X=(x_n)_{n\in \mathbb N}$ and $Y=(y_n)_{n\in \mathbb N}$.
To begin with, we will see that
${\textsl {Ind}}_{r,\textsf F}(X+Y)={\textsl {Ind}}_{r,\textsf F}(X)+{\textsl {Ind}}_{r,\textsf F}(Y).$
Since
\begin{align*}
(\Phi_{\rm ab})^{g(q)}_{\textsf F}\circ {\textsl{Ind}}_{r,\textsf F}(X+Y)
&={\textsf {Ind}}_{r,\textsf F}\circ (\Phi_{\rm ab})^{g(q)}_{\textsf  F}(X+Y)\\
&={\textsf {Ind}}_{r,\textsf F}\circ (\Phi_{\rm ab})^{g(q)}_{\textsf F}(X)+{\textsf {Ind}}_{r,\textsf F}\circ (\Phi_{\rm ab})^{g(q)}_{\textsf F}(Y)\\
&=(\Phi_{\rm ab})^{g(q)}_{\textsf F}\circ {\textsl{Ind}}_{r,\textsf F}(X)+(\Phi_{\rm ab})^{g(q)}_{\textsf F}\circ {\textsl{Ind}}_{r,\textsf F}(Y),
\end{align*}
our assertion follows from the injectivity of $(\Phi_{\rm ab})^{g(q)}_{\textsf F}$.
For any elements ${a}=(a_n)_{n\in \mathbb N}$ and ${b}=(b_n)_{n\in \mathbb N}$ of ${{\overline {W}}}^{g(q)}(A)$,
let $\pi: \textsf F\to A$ be the unique $\mathbb Z[q]$-algebra homomorphism with $x_n\mapsto a_n, y_n \mapsto b_n$
for all $n\in \mathbb N$.
Then we have
\begin{align*}
{\textsl{Ind}}_{r,A}((a_n)_{n\in \mathbb N}+(b_n)_{n\in \mathbb N})
&={\textsl{Ind}}_{r,A}((a_n)_{n\in \mathbb N}+(b_n)_{n\in \mathbb N})\\
&={\textsl{Ind}}_{r,A}({{\overline {W}}}^{g(q)}(\pi)(X)+{{\overline {W}}}^{g(q)}(\pi)(Y))\\
&={\textsl{Ind}}_{r,A}({{\overline {W}}}^{g(q)}(\pi)(X+Y)).
\end{align*}
Since
${\textsl{Ind}}_{r,A} \circ {{\overline {W}}}^{g(q)}(\pi)={{\overline {W}}}^{g(q)}(\pi)\circ {\textsl{Ind}}_{r,\textsf F},$
we finally have
\begin{align*}
{\textsl{Ind}}_{r,A}((a_n)_{n\in \mathbb N}+(b_n)_{n\in \mathbb N})
&={{\overline {W}}}^{g(q)}(\pi)\circ {\textsl{Ind}}_{r,\textsf F}(X)+{{\overline {W}}}^{g(q)}(\pi)\circ {\textsl{Ind}}_{r,\textsf F}(Y)\\
&={\textsl{Ind}}_{r,A}((a_n)_{n\in \mathbb N})+{\textsl{Ind}}_{r,A}((b_n)_{n\in \mathbb N})\\
&={\textsl{Ind}}_{r,A}((a_n)_{n\in \mathbb N})+{\textsl{Ind}}_{r,A}((b_n)_{n\in \mathbb N}).
\end{align*}
Moreover, it is easy to see that ${\textsl{Ind}}_{r}$ is a natural transformation.

{\bf Step 3:}
For the uniqueness, let $g(q)\ne 1$ and let $V_r$ be a natural transformation with the desired property.
Let $\textsf E=\mathbb Z[q][x_n: n\ge 1]$,
Let $X=(x_n)_{n\in \mathbb N}$ and $\textsf V=(\textsf v_n)_{n\in \mathbb N}:=V_{r,\textsf E}(X).$
From $\Phi_\textsf E^{g(q)}\circ V_{r,\textsf E}=V_{r,\textsf E}\circ \Phi_\textsf E^{g(q)}$ it follows that
$$\sum_{d|n}d(1-g(q)^{\frac nd})\textsf v_d^{\frac nd}
=r\sum_{e|\frac nr}e(1-g(q)^{\frac {n}{er}})x_e^{{\frac {n}{er}}}
=\sum_{d|n}d(1-g(q)^{\frac nd})x_{\frac dr}^{\frac nd}.$$
By the injectivity of $\Phi_\textsf E^{g(q)}$ we get $\textsf v_n=x_{\frac nr}$ for all $n\ge 1$.
Now, for any $\mathbb Z[q]$-algebra $A$ and any element $(a_n)_{n\in \mathbb N}\in {\overline {W}}^{g(q)}(A)$,
let $\hat \pi:\textsf E \to A$ be the morphism  uniquely determined by $x_n \mapsto a_n,\, (n\ge 1).$
From the commutativity of the diagram
\begin{equation*}
\begin{CD}
{\overline {W}}^{g(q)}(\textsf E) @>V_{r,\textsf E}>>{\overline {W}}^{g(q)}(\textsf E)\\
@V{{\overline {W}}^{g(q)}(\hat \pi )}VV  @VV{{\overline {W}}^{g(q)}(\hat \pi )}V \\
 {\overline {W}}^{g(q)}(A)@>>V_{r,A}> {\overline {W}}^{g(q)}(A)
\end{CD}
\end{equation*}
it follows that
$${V}_{r,A}((a_n)_{n\in \mathbb N})={V}_{r,A}( {\overline {W}}^{g(q)}(\hat \pi )(X))
= {\overline {W}}^{g(q)}(\hat \pi )\circ {V}_{r,\textsf E}(X)=
(a_{\frac nr})_{n\in \mathbb N}={\textsl {Ind}}_{r,A}((a_n)_{n\in \mathbb N}), $$
which shows that ${V}_{r,A}={\textsl {Ind}}_{r,A}$.
\qed
\vskip 2mm

In case of ${\overline \Lambda}^{g(q)}$ and $\Lambda^m\,(m\in \mathbb Z)$,
the $r$th inductions $V^{{\overline \Lambda}^{g(q)}}_r$ and $V^{\Lambda^m}_r\,(m\in \mathbb Z)$ are defined in such a way that
$V^{{\overline \Lambda}^{g(q)}}_{r,A}=\Theta_A\circ \textsl{Ind}_{r,A} \circ \Theta_A^{-1}$ and
$V^{\Lambda^m}_{r,A}=\Theta_A\circ \textsl{Ind}_{r,A} \circ \Theta_A^{-1}.$
In fact, it is easy to see that, for any object $A$,
\begin{align*}
&V^{{\overline \Lambda}^{g(q)}}_{r,A}(1+a_1t+a_2t^2+\cdots)=1+a_1t^r+a_2t^{2r}+\cdots, \\
&V^{\Lambda^m}_{r,A}(1+a_1t+a_2t^2+\cdots)=1+a_1t^r+a_2t^{2r}+\cdots.
\end{align*}

\subsection{Restriction}
Assume that ${\mathcal F}$ denotes one of the following functors
$$\textsf{gh}\circ F^{\textsf{{$\mathbb Z[q]$-Algebras}}}_{\textsf{Rings}},\,\,
\textsf{gh}\circ F^{{\textsf{$\mathbb Z[q]\otimes {\Psi}$-Rings}}}_{\textsf{Rings}},\,\,\text{ and }\textsf{gh}.$$
For each $r\ge 1$ and each object $A$, let
${\textsf {Res}}_{r,A}:\mathcal F(A)\to \mathcal F(A)$ be the map sending $(a_n)_{n\in \mathbb N}$ to
$(a_{rn})_{n\in \mathbb N}.$
Obviously ${\textsf {Res}}_{r,A}$ is a ring homomorphism.
Moreover, it is easy to show that the assignment $A\mapsto {\textsf {Res}}_{r,A}$ induces a natural transformation
${\textsf {Res}}_{r}:\mathcal F \to \mathcal F$.
\subsubsection{}
Here we introduce restrictions on ${\overline {W}}^{g(q)}$ and ${W}^m,\, (m\in \mathbb Z)$.
\begin{prop}\label{uniquedness of restriction q-witt }
Given a polynomial $g(q)\in \mathbb Z[q]$, there exists a natural transformation
${\textsl {Res}}_{r}^{{\overline {W}}^{g(q)}}:{\overline {W}}^{g(q)}\to {\overline {W}}^{g(q)}$
with
${\overline \Phi}^{g(q)}\circ {\textsl {Res}}_{r}^{{\overline {W}}^{g(q)}}={\textsf {Res}}_{r}\circ {{\overline \Phi}^{g(q)}}.$
Furthermore, if $g(q)\ne 1$, then ${\textsl {Res}}_{r}^{{\overline {W}}^{g(q)}}$ is uniquely determined
by the above commutativity.
\end{prop}

\pf
Let $\textsf F=\mathbb Z[q][x_n:n\ge 1]$.
Let $X=(x_n)_{n\in \mathbb N}$ and define $\textsf R_n, (n\ge 1)$
via the following recursive relations:
\begin{equation}\label{how to define restriction for w(g(q))}
n\textsf R_n
=\sum_{e| nr}e\left[\frac {nr}{e}\right]_{g(q)}x_e^{{\frac {nr}{e}}}
-\sum_{d|n \atop d<n}d\left[\frac nd\right]_{g(q)}\textsf R_d^{\frac nd}, \quad (n\ge 1).
\end{equation}
Obviously $\textsf R_n \in \mathbb Q[g(q)][x_d: d|nr]$.
However, we claim that $\textsf R_n \in \mathbb Z[g(q)][x_d: d|nr]$ for all $n\ge 1$.

First, assume that $g(q)=1$:
Eq.\eqref{how to define restriction for w(g(q))} is simplified as
$$\sum_{d|n}\textsf R_d^{\frac nd}
=r\sum_{e| nr}x_e^{{\frac {nr}{e}}}.$$
Since $\textsf R_1=r\sum_{e| r}x_e^{{\frac {r}{e}}},$
our assertion holds when $n=1$.
Now the claim can be verified by using the mathematical induction on $n$.

Next, let us consider the case $g(q)\ne 1$.
Since $$\textsf R_1=\sum_{e| r}e\left[\frac {r}{e}\right]_{g(q)}x_e^{{\frac {r}{e}}},$$
our claim holds when $n=1$.
If $n>1$, assume that the claim holds for all positive integers less than $n$.
Then it is easy to see that the right hand side of Eq.\eqref{how to define restriction for w(g(q))}
is in $\mathbb Z[g(q)][x_d: d|nr]$.
Thus, for the verification of our claim, we have only to see that the right hand side of Eq.\eqref{how to define restriction for w(g(q))}
is divisible by $n$. To do this, we will see that the right hand side of Eq.\eqref{how to define restriction for w(g(q))} multiplied by $1-g(q)$
is divisible by $n(1-g(q))$.
But it can be done exactly in the same way as in the proof of Lemma \ref{key lemma in the proof of existence of q deformation of Witt rings}.
The only nontrivial part is to see that $n(1-{g(q)})\textsf R_n$ is divisible by $n$.
Indeed this is true since, for every prime divisor $p$ of $n$, we have
\begin{align*}
n(1-{g(q)})\textsf R_n
&=\sum_{e| nr}e(1-{g(q)}^{\frac {nr}{e}})x_e^{{\frac {nr}{e}}}-\sum_{d|n \atop d<n}d(1-{g(q)}^{\frac nd})\textsf R_d^{\frac nd} \\
& \equiv \sum_{e| \frac {nr}{p}}e(1-({g(q)}^p)^{\frac {nr}{pe}})(x_e^p)^{{\frac {nr}{pe}}}
-\sum_{d|\frac np}d(1-{g(q)}^{\frac {n}{d}})\textsf R_d^{\frac nd} \\
&\equiv \sum_{d|\frac np}d(1-{g(q)}^{\frac {n}{d}})
((\textsf R_d|_{g(q) \mapsto g(q)^p \atop x_e \mapsto x_e^p})^{\frac{n}{pd}}-(\textsf R_d^p)^{\frac {n}{pd}})
\\
&\equiv 0 \pmod{p^{\nu_p(n)}} \quad \text{(by Lemma \ref{key modular equivalence 1}(a))}.
\end{align*}
This shows that $\textsf{R}_n \in \mathbb Z[g(q)][x_d: d|n, d<n ]$.
For any $\mathbb Z[q]$-algebra $A$, consider the map
$${\textsl {Res}}_{r,A}^{{\overline {W}}^{g(q)}}:{\overline {W}}^{g(q)}(A)\to {\overline {W}}^{g(q)}(A),
\quad (a_n)_{n \in \mathbb N}\mapsto (\textsf R_n|_{x_d \mapsto a_d})_{n \in \mathbb N}.$$
Mimicking the proof of Proposition \ref{uniquness of induction B},
one can see that the assignment $A\mapsto {\textsl {Res}}_{r,A}^{{\overline {W}}^{g(q)}}$ induces a natural transformation with the desired property.
\qed

\vskip 2mm
Let $g(q)=q$.
Then Proposition \ref{uniquedness of restriction q-witt } produces the following corollary
when $q$ is specialized into $m$.

\begin{cor}\label{restriction of specialized Witt }
Let $m$ be an integer. For each positive integer $r$, there exists a natural transformation
${\textsl {Res}}_{r}^{W^m}:W^m\to W^m$ with $\Phi^m \circ {\textsl {Res}}_{r}^{W^m}={\textsf {Res}}_{r}\circ \Phi^m$.
Furthermore, if $g(q)\ne 1$, then ${\textsl {Res}}_{r}^{W^m}$ satisfying this commutativity is uniquely determined.
\end{cor}

\begin{rem}{\rm
In case of ${\overline \Lambda}^{g(q)}$ and $\Lambda^m,(m\in \mathbb Z)$,
the $r$th restrictions $F^{{\overline \Lambda}^{g(q)}}_r$ and $F^{\Lambda^m}_r,\,(m\in \mathbb Z),$ are defined in such a way that
$F^{{\overline \Lambda}^{g(q)}}_{r,A}=\Theta_A\circ \textsl{Res}_{r,A} \circ \Theta_A^{-1}$ and
$F^{\Lambda^m}_{r,A}=\Theta_A\circ \textsl{Res}_{r,A} \circ \Theta_A^{-1}.$
It is easy to see that, for any $\mathbb Z[q]$-algebra $A$,
\begin{align*}
&F^{{\overline \Lambda}^{g(q)}}_{r,A}\left(\prod_{i=1}^{\infty}\frac{1}{1-a_nt^n}\right)=\prod_{i=1}^{\infty}\frac{1}{1-\textsf R_n|_{x_d\mapsto a_d}t^n}.
\end{align*}
}\end{rem}

\subsubsection{}
One can show that statements analogous to
Proposition \ref{uniquedness of restriction q-witt } and Corollary \ref{restriction of specialized Witt }
hold for ${\overline {\textsf B}}^{g(q)}$ and $\textsf B^m$, respectively.

\begin{prop}\label{uniquedness of restriction q-necklace func}
{\rm (a)}
Given a polynomial $g(q)\in \mathbb Z[q]$, there exists a natural transformation
$\textsl {Res}_r^{B^{g(q)}}:{\overline {\textsf B}}^{g(q)}\to {\overline {\textsf B}}^{g(q)}$
with
${\overline \varphi}^{g(q)}\circ {\textsl {Res}}_{r}^{{\overline {\textsf B}}^{g(q)}}={\textsf {Res}}_{r}\circ {{\overline \varphi}^{g(q)}}.$
Furthermore, if $g(q)\ne 1$, then $\textsl {Res}_r^{B^{g(q)}}$ is uniquely determined
by the above commutativity.

{\rm (b)}
For each integer $m$, there exists a natural transformation
${\textsl {Res}}_{r}^{\textsf B^m}:\textsf B^m \to \textsf B^m$ with $\varphi^m \circ {\textsl {Res}}_{r}^{\textsf B^m}={\textsf{Res}}_{r}\circ \varphi^m.$
Furthermore, if $m\ne 1$, then ${\textsl {Res}}_{r}^{\textsf B^m}$ satisfying this commutativity is uniquely determined.

{\rm (c)}
For each $\Psi$-ring $A$ and $(a_n)_{n \in \mathbb N}\in \textsf B^m(A)$, we have
\begin{equation}\label{closed form of restriction3 }
{\textsl {Res}}_{r}^{\textsf B^m}((a_n)_{n \in \mathbb N})
=\left(\frac {1}{n}\sum_{e|nr}
\left(e\sum_{d|(n,\frac{nr}{e})}\mu_m(d)\left[\frac{nr}{de}\right]_m \right)
\Psi^{\frac{nr}{e}}(a_e)\right)_{n \in \mathbb N}.
\end{equation}
\end{prop}

\pf
The first two statements can be proved by slightly modifying the proof of Proposition\ref{uniquedness of restriction q-witt }.
So we here prove only (c).
Note that
$$\sum_{d|n}d\left[\frac nd\right]_{m}\Psi^{\frac nd}(\textsf R_d)
=\sum_{e| nr}e\left[\frac {nr}{e}\right]_{m}\Psi^{{\frac {nr}{e}}}(x_e), \quad (n\ge 1).$$
In view of Lemma \ref{q mobius ftn psi-operation}, one has
\begin{align*}
{{\textsf R}}_n
=\frac 1n \sum_{d|n}\mu_m(d)\Psi^d \left(\sum_{e| \frac{nr}{d}}e\left[\frac {nr}{ed}\right]_{m}
\Psi^{{\frac {nr}{ed}}}(x_e)\right)
=\frac 1n\sum_{e|nr}e \left(\sum_{d|(n,\frac{nr}{e})}\mu_m(d)\left[\frac {nr}{ed}\right]_{m}\right)
\Psi^{\frac{nr}{e}}(x_e).
\end{align*}
As a consequence,
\begin{equation}\label{closed fomular for restructuin B(m) case}
{{\textsf R}}_n|_{x_d \mapsto a_d}
=\frac 1n\sum_{e|nr}e\left(\sum_{d|(n,\frac{nr}{e})}\mu_m(d)\left[\frac {nr}{ed}\right]_{m}\right)\Psi^{\frac{nr}{e}}(a_e),
\end{equation}
as required.
\qed

\begin{rem}{\rm
If $m=0$, then Eq.\eqref{closed fomular for restructuin B(m) case} reduces
$${{\textsf R}}_n|_{x_d \mapsto a_d}
=\frac 1n\sum_{[r,e]=nr}e\Psi^{\frac{nr}{e}}(a_e).$$
}\end{rem}

\subsection{Unitalness theorems}
Given a ring-valued functor, say ${\mathcal F}:\mathcal C\to \mathcal D$,
one of the most fundamental problems may be to determine whether ${\mathcal F}(A)$ is unital or not
for each object $A$ in $\mathcal C$.
In the rest of this section, we deal with this problem for our generalized functors.

\begin{prop}\label{Unity existence theorem}
Suppose $g(q)\ne 1.$ Then the following hold.

{\rm (a)}
Let $A$ be a $\mathbb Z[q]$-algebra. Then
${\overline {W}}^{g(q)}(A)$ is unital if and only if $1-g(q)$ is a unit in $A$.

{\rm (b)}
Let $A$ be an object of the category {\textsf{$\mathbb Z[q]\otimes {\Psi}$-Rings}}. Then
${{{\overline {\textsf B}}}}^{g(q)}(A)$ is unital if and only if $1-g(q)$ is a unit in $A$.
\end{prop}

\pf
(a) First, assume that ${\overline {W}}^{g(q)}(A)$ is unital.
Let $(e_n)_{n\in \mathbb N}$ be the unity of ${\overline {W}}^{g(q)}(A)$ and let
$(a_n)_{n\in \mathbb N}$ be an arbitrary element in ${\overline {W}}^{g(q)}(A)$.
Since $(e_n)_{n\in \mathbb N}\cdot (a_n)_{n\in \mathbb N}=(a_n)_{n\in \mathbb N}$, it holds that $(1-g(q))a_1e_1=a_1$.
In particular, by letting $a_1=1$, we have $(1-g(q))e_1=1$
So, $1-g(q)$ is a unit in $A$.
For the converse,
let us define $\textsf E_n \in \mathbb Q(q), (n\ge 1)$ via the following recursions:
\begin{equation}\label{existence of unity00}
\sum_{d|n}d (1-g(q)^{\frac{n}{d}})\textsf E_d^{\frac nd}=1,\,\, (n\ge 1).
\end{equation}
We claim that $(1-g(q))^{n}\textsf E_n \in \mathbb Z[q]$.
The case where $n=1$ is trivial since $(1-g(q))\textsf E_1=1$.
If $n>1$, assume $(1-g(q))^{i}\textsf E_i \in \mathbb Z[q]$ for all positive integers $i<n$.
Let
$$\textsf E_i:=\frac{{\hat {\textsf E}}_i}{(1-g(q))^i}$$
for $1\le i < n.$
Multiplying both sides of Eq.\eqref{existence of unity00} by $(1-g(q))^{n}$ yields that
\begin{align}\label{condition for unity existence witt ring}
n(1-g(q))^{n+1}\textsf E_n
=(1-g(q))^{n}-\left(\sum_{d|n \atop d<n}d (1-g(q)^{\frac{n}{d}}){\hat {\textsf E}}_d^{\frac nd}\right).
\end{align}
Since ${\hat {\textsf E}}_d\in \mathbb Z[q]$, the right hand side is contained in $(1-g(q))\mathbb Z[q]$ and thus
$n(1-g(q))^{n+1}\textsf E_n$ is divisible by $1-g(q).$
To see that $n(1-g(q))^{n+1}\textsf E_n$ is divisible by $n$,
choose any prime divisor $p$ of $n$.
Notice
\begin{align*}
\sum_{d|n \atop d<n}d (1-g(q)^{\frac{n}{d}}){\hat {\textsf E}}_d^{\frac nd}
&\equiv \sum_{d|\frac np}d (1-(g(q)^p)^{\frac{n}{dp}})({\hat {\textsf E}}_d^p)^{\frac {n}{dp}} \\
&\equiv \sum_{d|\frac np}d (1-g(q^p)^{\frac{n}{dp}})({\hat {\textsf E}}_d|_{q \mapsto q^p})^{\frac {n}{dp}} \pmod {p^{\nu_p(n)}},
\end{align*}
which holds since
$$\sum_{d|\frac np}d (1-g(q^p)^{\frac{n}{dp}})(({\hat {\textsf E}}_d|_{q \mapsto q^p})^{\frac {n}{dp}}-({\hat {\textsf E}}_d^p)^{\frac {n}{dp}})
\equiv 0 \pmod {p^{\nu_p(n)}}.$$
On the other hand, in view of Eq.\eqref{existence of unity00}, one has that
$$\sum_{d|\frac np}d (1-g(q^p)^{\frac{n}{dp}})({\hat {\textsf E}}_d|_{q \mapsto q^p})^{\frac {n}{dp}}=(1-g(q^p))^{\frac np}.$$
As a consequence,
$$n(1-g(q))^{n+1}\textsf E_n\equiv (1-g(q))^n-(1-g(q^p))^{\frac np} \equiv 0 \pmod {p^{\nu_p(n)}}.$$
The last equivalence follows from Lemma \ref{key modular equivalence 1}(a).
Now, we are ready to prove $(1-g(q))^{n}\textsf E_n\in \mathbb Z[q].$
Let $p'$ be a prime divisor of $n(1-g(q))$.
In case where $(p',n)=1$ or $(p',1-g(q))=1$, our assertion is trivial since
$n(1-g(q))((1-g(q))^{n}\textsf E_n )$ is divisible by $n$ and $1-g(q)$.
So we assume $(p',n)=(p',1-g(q))=p'$.
For every divisor $d$ of $n$, it holds that
$$d(1-g(q)^{\frac{n}{d}})\equiv {p'}^{\nu_p(n(1-g(q))}$$
(refer to the proof of Lemma \ref{key lemma in the proof of existence of q deformation of Witt rings}).
It says that the right hand side of Eq.\eqref{condition for unity existence witt ring}
is divisible by ${p'}^{\nu_p(n(1-g(q))}$.
So our claim is verified.
Finally, set $e_n=\textsf E_n \cdot 1_A$ for all $n\ge 1$.
By construction it is obvious that
$(e_n)_{n\in \mathbb N}$ is the unity of ${\overline {W}}^{g(q)}(A).$

(b) The ``only if" part can be shown as in the same way as in (a).
For the converse, assume that $1-g(q)$ is a unit.
For all $n\ge 1$ define $E_n \in \mathbb Q(q)$ to be
$$\begin{cases}
(1-g(q))^{-1} & \text{ if }n=1,\\
0 & \text{ otherwise.}
\end{cases}$$
Set $e_n=\textsf E_n \cdot 1_A$ for all $n\ge 1$.
Then
$$\sum_{d|n}d (1-g(q)^{\frac nd})\Psi^{\frac nd}(e_d)=1, \quad (n\ge 1), $$
which means that
$(e_n)_{n\in \mathbb N}$ is the unity of ${{\overline {\textsf B}}}^{g(q)}(A).$
\qed

\vskip 2mm
The following corollary can be deduced from Proposition \ref{Unity existence theorem} by specializing $q$ into $m$.
\begin{cor}
Suppose that $m$ is an integer not equal to $1$. Then the following hold.

{\rm (a)}
Let $A$ be an commutative ring with unity. Then
$W^m(A)$ is unital if and only if $1-m$ is a unit in $A$.

{\rm (b)}
Let $A$ be a $\Psi$-ring. Then
$\textsf B^m(A)$ is unital if and only if $1-m$ is a unit in $A$.
\end{cor}

\section{Relationship between Witt rings and necklace rings}
\label{Relationship between Witt rings and necklace rings}

\subsection{Connection between ${\overline {W}}^{g(q)}$ and ${{\overline {\textsf B}}}^{g(q)}$}
Let {\textsf{$\mathbb Z[q]\otimes {\lambda}$-Rings}} be the subcategory of
{\textsf{$\mathbb Z[q]\otimes {\Psi}$-Rings}}
such that the objects are $\mathbb Z[q]$-algebras equipped with the $\ld$-ring structure with $\psi^n=\Psi^n$ for all $n\ge 1$
and the morphisms are $\mathbb Z[q]$-algebra homomorphisms preserving $\ld$-operations.
Therefore,
$$\psi^n(g(q))=g(q^n),  \quad g(q)\in \mathbb Z[q], \,\, n\ge 1.$$
Let $\textsf E$ be the object of {\textsf{$\mathbb Z[q]\otimes {\ld}$-Rings}}
freely generated by $x$, thus
$$\textsf E=\mathbb Z[q][\ld^n(x): n\ge 1].$$
For each positive integer $n$ define ${\textsf M}^{g(q)}(x,n)$ via the following recursive relations:
$$\sum_{d|n}d(1-g(q)^{\frac nd})\Psi^{\frac nd}({\textsf M}^{g(q)}(x,d))=(1-g(q)^n)x^n,\quad (n\ge 1).$$
Let $p$ be a prime divisor of $n$.
Then, in view of \eqref{key modular relation}, we have
\begin{align*}
n(1-g(q))\textsf M^{g(q)}(x,n)
&=(1-g(q)^n)x^n-\sum_{d|n}d(1-g(q)^{\frac nd})\Psi^{\frac nd}({\textsf M}^{g(q)}(x,d))\\
&\equiv(1-g(q)^n)x^n-\sum_{d|\frac np}d(1-g(q^p)^{\frac {n}{dp}})\Psi^{\frac nd}({\textsf M}^{g(q)}(x,d))\pmod{p^{\nu_p(n)}}.
\end{align*}
Applying the following congruence
\begin{align*}
\sum_{d|\frac np}d(1-g(q^p)^{\frac {n}{dp}})\Psi^{\frac nd}({\textsf M}^{g(q)}(x,d))
&\equiv \Psi^p(\sum_{d|\frac np}d(1-g(q)^{\frac {n}{dp}})\Psi^{\frac {n}{dp}}({\textsf M}^{g(q)}(x,d)))\\
&\equiv \Psi^p((1-g(q)^{\frac np})x^{\frac np})\\
&\equiv (1-g(q^p)^{\frac np})\Psi^p(x^{\frac np}) \pmod{p^{\nu_p(n)}}
\end{align*}
yields that
$$n(1-g(q))\textsf M^{g(q)}(x,n)\equiv(1-g(q)^n)(x^n-\Psi^p(x^{\frac np})) \pmod{p^{\nu_p(n)}}.$$
Let $n=p^{\nu_p(n)}n'$. Then Lemma \ref{to define Teichmuller map} implies that
$$M(x^{n'},p^{\nu_p(n)})=\frac{1}{p^{\nu_p(n)}}(x^n-\psi^p(x^{\frac np})) \in \textsf E,$$
and thus
$$x^n-\psi^p(x^{\frac np})\equiv 0 \pmod{p^{\nu_p(n)}}.$$
Let $R$ be an arbitrary object of the category  {\textsf{$\mathbb Z[q]\otimes {\ld}$-Rings}}.
Given an element $r\in R$ let $\pi: \textsf E \to R$ be the unique $\mathbb Z[q]$-algebra homomorphism preserving $\ld$-operations
determined by $x\mapsto r$.
Define ${\textsf M}^{g(q)}(r,n)$ to be $\pi({\textsf M}^{g(q)}(x,n))$ for all $n\ge 1$
and consider the map
$${\textsf M}^{g(q)}: R \to {{\overline {\textsf B}}}^{g(q)}(R), \quad r \mapsto ({\textsf M}^{g(q)}(r,n))_{n\in \mathbb N}.$$
Finally, let ${\overline \tau}^{g(q)}_R: {\overline {W}}^{g(q)}(R) \to {{\overline {\textsf B}}}^{g(q)}(R)$ be the map defined by
$$(a_n)_{n\in \mathbb N}\mapsto \sum_{n\ge 1}{\textsl {Ind}}_{n,R}({\textsf M}^{g(q)}(a_n)).$$
Then
$${\overline \tau}^{g(q)}_R
((a_n)_{n\in \mathbb N})=\left(\sum_{d|n}{\textsf M}^{g(q)}(a_d, \frac nd )\right)_{n\in \mathbb N}.$$

\begin{thm} \label{relation of W(g(q)) and B(g(q))}
The assignment $R\mapsto {\overline \tau}^{g(q)}_R$ induces a natural isomorphism
$${\overline \tau}^{g(q)}: {\overline {W}}^{g(q)} \circ F^{{\textsf{$\mathbb Z[q]\otimes {\ld}$-Rings}}}_{\textsf{Rings}} \to
{{\overline {\textsf B}}}^{g(q)}|_{{\textsf{$\mathbb Z[q]\otimes {\ld}$-Rings}}}.$$
\end{thm}

\pf
To begin with, we show that, for every object $R$ in {\textsf{$\mathbb Z[q]\otimes {\ld}$-Rings}},
${\overline \tau}^{g(q)}_R$ is a ring isomorphism with ${\overline \Phi}^{g(q)}_R={\overline \varphi}^{g(q)}_R\circ {\overline \tau}^{g(q)}_R.$
Note that
$${\overline \tau}^{g(q)}_R((a_n)_{n\in \mathbb N})
=\left(a_n+\sum_{d|n \atop d<n}{\textsf M}^{g(q)}(a_d, \frac nd )\right)_{n\in \mathbb N}.$$
Since the coefficient of $a_n$ is $1$, one can see in a recursive way that ${\overline \tau}^{g(q)}_R$ is bijective.
To see that ${\overline \Phi}^{g(q)}_R={\overline \varphi}^{g(q)}_R\circ {\overline \tau}^{g(q)}_R,$
choose any element $(a_n)_{n\in \mathbb N}\in {\overline {W}}^{g(q)}(R)$.
Then we have
\begin{align*}
{\overline \varphi}^{g(q)}_R\circ \sum_{r\ge 1}{\textsl{Ind}_r}(M^{g(q)}(a_r))
&\stackrel{\text{by Prop. \ref{uniquness of induction B}}}{=}\sum_{r\ge 1}{\textsl{Ind}_r}\circ {\overline \varphi}^{g(q)}_R(M^{g(q)}(a_r))\\
&\quad \,\,\,\,\,=\sum_{r\ge 1}{\textsl{Ind}_r}(((1-g(q)^{n})a_r^n)_{n\ge 1}).
\end{align*}
Note that
$${\textsl{Ind}_r}(((1-g(q)^{n})a_r^n)_{n\ge 1})
=(b_{r,n})_{n \in \mathbb N},$$
where
$$
b_{r,m}=
\begin{cases}
r(1-g(q)^{\frac nr})a_r^{\frac nr} & \text{ if } r|n, \\
0 & \text{ otherwise.}
\end{cases}$$
Consequently the $n$th component of $\sum_{r\ge 1}{\textsl{Ind}_r}(M^{g(q)}(a_r))$ is given by
$\sum_{d|n}d(1-g(q)^{\frac nd})a_d^{\frac nd}.$
Finally we should verify the commutativity condition.
It is straightforward since, for any $\mathbb Z[q]$-algebra homomorphism $f:R\to S$ preserving $\ld$-operations,
$$\textsf M^{g(q)}(f(a), n)=f(\textsf M^{g(q)}(a, n)), \quad a\in R, n\ge 1.$$
So we are done.
\qed

\vskip 2mm
Let $R$ be an object of the category {\textsf{$\mathbb Z[q]\otimes {\ld}$-Rings}}. Then the diagram
\begin{equation*}
\begin{CD}
{{{\overline {\textsf B}}}}^{g(q)}|_{{\textsf{$\mathbb Z[q]\otimes {\ld}$-Rings}}}(R) @<{\overline \tau}^{g(q)}_R<<
{\overline {W}}^{g(q)}(R)@>\Theta_R>>
\Lambda^0 (R)\\
@V{{\overline \varphi}^{g(q)}_R}VV  @V{{\overline \Phi}^{g(q)}_R}VV  @V{{\overline \Upsilon}^{g(q)}_R}VV \\
(\textsf{gh}(R) @<{\rm id}<<
(\textsf{gh}(R) @>{\it int}_R>>
(\textsf{Gh}(R)
\end{CD}
\end{equation*}
is commutative. From the commutativity it follows that
\begin{equation}\label{symmetric map}
\Theta_R \circ ({\overline \tau}^{g(q)}_R)^{-1}((a_n)_{n\in \mathbb N})
=\prod_{n\ge 1}\exp \left(\sum_{r\ge 1}\frac 1r(1-g(q)^r)\psi^{r}(a_n)t^{nr}\right).
\end{equation}
In particular, if $\psi^{r}={\rm id}$ for all $r\ge 1$, then Eq.\eqref{symmetric map} is reduced to
\begin{equation}\label{symmetric map1}
\Theta_R \circ ({\overline \tau}^{g(q)}_R)^{-1}((a_n)_{n\in \mathbb N})
=\prod_{n\ge 1}\left(\frac{1-g(q)t^n}{1-t^n}\right)^{a_n}.
\end{equation}

\begin{rem}\label{Henry Kim}{\rm
Recently Kim and Lee \cite{KL} introduced the product identities of the following form:
\begin{align}\label{Henry KIM and Kyuhwan Lee}
\prod_{n\ge 1}\left(\frac{1-t^n}{1-qt^n}\right)^{a_n}=1+\sum_{n\ge 1}b_nt^n,
\end{align}
where $a_n\in \mathbb Z$ and $b_n\in \mathbb Z[q]$ for all $n\ge 1$.
In order to interpret them ring-theoretically, let
$R=\mathbb Z[q]$ be the $\ld$-ring with $\psi^n(q)=q^n,\,(n\ge 1).$
Let $$\bar \Lambda^q(R)=\{1+\sum_{n\ge 1}a_nt^n: a_n\in R\}$$
be the ring such that
the ring operations are defined in such a way that
\begin{align*}
&f(t)+g(t):=\rho(\rho^{-1}(f(t))+\rho^{-1}(g(t))), \\
&f(t)*g(t):=\rho(\rho^{-1}(f(t))*\rho^{-1}(g(t))),
\end{align*}
where
$$\rho:\Lambda^q(R)\to \bar \Lambda^q(R), \,\,\,f(t)\mapsto \frac{1}{f(t)}.$$
In view of Eq.\eqref{symmetric map1}, if $a_n\in \mathbb Z$ for all $n\ge 1$, then
$$\rho\circ \Theta_R \circ ({\overline \tau}_R^q)^{-1}((a_n)_{n\in \mathbb N})
=\prod_{n\ge 1}\left(\frac{1-t^n}{1-qt^n}\right)^{a_n}.$$
Finally, let
$$1+\sum_{n\ge 1}b_nt^n:=\rho\circ \Theta_R \circ ({\overline \tau}_R^q)^{-1}((a_n)_{n\in \mathbb N}).$$
Since $\rho, \Theta_R, {\overline \tau}_R^q$ are isomorphisms, $b_n\in \mathbb Z[q]$ for all $n\ge 1$ and thus we obtain
Eq.\eqref{Henry KIM and Kyuhwan Lee}.
}\end{rem}

\subsection{Connection between $W^m$ and $\textsf B^m$}
Let us define an $\mathbb N\times \mathbb N$ matrix
$\zeta_m$ with entries in $\mathbb Z$ by
$$\zeta_m(a,b)=
\begin{cases}
\left[\frac ab \right]_m
& \text{ if }b|a, \\
0 & \text { otherwise.}
\end{cases}$$
It is a lower-triangular matrix
with the diagonal elements $1$.
Denote by $\mu_m$ the inverse of $\zeta_m$.
As in Section \ref{Structure constants},
for any chain
$C=(a=x_r \prec \cdots \prec x_0=b)\in \mathcal C(a,b),$ we let
$${\rm wt}_{m}(C):=(-1)^r\left[\frac{x_0}{x_1}\right]_{m}\cdots \left[\frac{x_{r-1}}{x_r}\right]_{m}.$$
Then it is easy to see
$$\mu_m(b,a)=\sum_{C\in \mathcal C(a,b)}{\rm wt}_m(C)$$
and $\mu_m(b,a)=\mu_m(b',a')$ if $b/a={b'}/{a'}$.
For each positive integer $n$, let us define $\mu_m(n)$ by $\mu_m(n,1)$.
Then the entries of $\mu_m$ have the following properties.
To begin with, let us introduce some notations.
Given a positive integer $n$, a {\it multiplicative composition} of $n$, denoted by $\alpha\vDash n$, is an ordered sequence $(\alpha_1, \alpha_2, \ldots, \alpha_l)$
of positive integers with $\alpha_1\cdot \alpha_2 \cdots \cdot\alpha_l=n$ and $1\le l\le n$.
In this case, $l$ is called the {\it length of $\alpha$}.
Set
$$[\alpha]_m:=[\alpha_1]_m [\alpha_2]_m\cdots [\alpha_l]_m.$$

\begin{lem}\label{q mobius ftn}
Let $a$ be a positive integer.
For any multiple $n$ of $a$, we have

{\rm (a)} $\mu_m(n,a)=\sum_{\alpha\vDash \frac na}(-1)^{l(\alpha)-1}[\alpha]_m$.

{\rm (b)} $\mu_m(n,a)=\mu_m(ln,la)$ for every positive integer $l$.

\end{lem}
\pf
(a) Note that the $(n,a)$th entry of $\zeta_m \mu_m$ equals
$$\sum_{a|d|n}\left[\frac nd \right]_m \mu_m(d,a)
=\begin{cases}
1 & \text{ if }a=n,\\
0 & \text{ otherwise.}
\end{cases}$$
If $a=n$, then the desired result obviously holds.
Assume that the assertion holds for all the multiples of $a$ less than $n$.
Then
\begin{align*}
\mu_m(n,a)
&=-\sum_{a|d|n \atop d<n} \left[\frac nd\right]_m (\sum_{\alpha \vDash \frac da}(-1)^{l(\alpha)-1}[\alpha]_m)\\
&=\sum_{\beta \vDash \frac na}(-1)^{l(\beta)-1}[\beta]_m \quad \text{(by letting $\beta=(\alpha_1,\ldots, \alpha_{l(\alpha)}, n/d)$}.
\end{align*}
(b) This follows from (a).
\qed

\vskip 2mm
Define $\mu_m(n)$ to be $\mu_m(n,1),\, (n\ge 1)$.
From Lemma \ref{q mobius ftn} it follows that
$$\mu_m(n,1)=\sum_{\alpha\vDash n}(-1)^{l(\alpha)-1}[\alpha]_m=\mu_m(nl,l)$$
for all $l\ge 1$.
For $m\ne 1$, let us define another $\mathbb N\times \mathbb N$ matrix
$\tilde \zeta_m$ by
$$\tilde\zeta_m(a,b)=
\begin{cases}
b(1-m^{\frac ab})\Psi^{\frac nd}
& \text{ if }b|a, \\
0 & \text { otherwise.}
\end{cases}
$$
Let $\tilde \mu_m$ be the inverse of $\tilde\zeta_m$.
Then we can derive the following lemma.

\begin{lem}\label{q mobius ftn psi-operation}
Let $a$ be any positive integer.
For any multiple $n$ of $a$, we have that
$$\tilde\mu_m(n,a)=\frac{\mu_m(n,a)}{n(1-m)} \Psi^{\frac na}.$$
\end{lem}

\pf
The $(n,a)$th entry of $\tilde \zeta_m \tilde \mu_m$ equals
$$\sum_{a|d|n}d(1-m^{\frac nd})\Psi^{\frac nd} \circ \tilde \mu_m(d,a)
=\begin{cases}
{\rm id} & \text{ if }a=n,\\
0 & \text{ otherwise.}
\end{cases}$$
If $a=n$, then our assertion is trivial since
$$\tilde \mu_m(a,a)=\frac {1}{a(1-m)} {\rm id}.$$
Now, we assume the our assertion holds for all multiples of $a$ less than $n$.
Then
\begin{align*}
n(1-m)\tilde\mu_m(n,a)
&=-\sum_{a|d|n \atop d<n}d(1-m^{\frac nd})\Psi^{\frac nd} \circ \frac {1}{d(1-m)}
\left(\sum_{\alpha \vDash \frac da}(-1)^{l(\alpha)-1}[\alpha]_m \right)\Psi^{\frac da}\\
&=\left(\sum_{\beta \vDash \frac na}(-1)^{l(\beta)-1}[\beta]_m \right) \Psi^{\frac na} .
\end{align*}
Now the desired result follows from Lemma \ref{q mobius ftn}(a).
\qed

\vskip 2mm
Let $R$ be a torsion free $\Psi$-ring.
For $x\in R$, let
$\textsf{M}^m(x)$ be the transpose of
$$\tilde\mu_m \begin{pmatrix}\vdots\\(1-m^n)x^n \\ \vdots\end{pmatrix}_{n\in \mathbb N}.$$
Let $\textsf{M}^m(x,n)$ be the $n$th component of $\textsf{M}^m(x)$, that is, $\textsf{M}^m(x)=(\textsf{M}^m(x,n))_{n\in \mathbb N}$.
It should be noticed that $\textsf{M}^m(x,n)$ is contained in $\mathbb Q\otimes R$, not in $R$, in general.
For all $n \ge 1$, we have
\begin{equation}\label{How to define expoenntial map}
\begin{aligned}
\textsf{M}^m(x,n)
&=\sum_{d|n}\tilde \mu_m(n,d)((1-m^d)x^d)\\
&=\frac{1}{n(1-m)}\sum_{d|n}\mu_m(n,d)\Psi^{\frac nd}((1-m^d)x^d) \quad \text{ (by Lemma \ref{q mobius ftn psi-operation})}\\
&=\frac {1}{n}\sum_{d|n}\mu_m(d)\left[ \frac nd \right]_m\Psi^{d}(x^{\frac nd}).
\end{aligned}
\end{equation}
Utilizing the relation
$$\tilde \zeta_m \textsf{M}^m(x)^T =((1-m^n)x^n)^T_{n\in \mathbb N},$$
we can derive the following recursive formula for $\textsf{M}^m(x,n)$:
$$\sum_{d|n}d \left[\frac nd \right]_m\Psi^{\frac nd}(\textsf{M}^m(x,d))=[n]_mx^n,\,\, (n\ge 1).$$

Let $R$ be an arbitrary $\ld$-ring and let $m$ be an arbitrary integer.
For $r\in R$ and each positive integer $n$, let us define $\textsf{M}^m(r,n)$.
To do this, we first consider the free $\ld$-ring generated by $x$, denoted by $F(x)$.
It is easy to see that
$F(x)=\mathbb Z[\ld^k(x): k \ge 1 ]$
(for instance, see \cite{AT}).
It is well known that $F(x)$ is isomorphic to $\Lambda_\mathbb Z$, the ring of symmetric functions over $\mathbb Z$ in $x_1,x_2,\ldots$,
where the explicit isomorphism is given by
$$F(x)\to \Lambda, \, \ld^k(x) \mapsto e_k(x),\quad  (k\ge 1).$$
For all $n\ge 1$, let
\begin{equation}\label{How to define expoenntial map 1}
\textsf{M}^m(x,n):=\frac {1}{n}\sum_{d|n}\mu_m(d)\left[ \frac nd \right]_m \psi^{d}(x^{\frac nd}).
\end{equation}
By definition it is straightforward that
$$\sum_{d|n}d \left[\frac nd \right]_m\psi^{\frac nd}(\textsf{M}^m(x,d))=[n]_mx^n.$$

\begin{lem}\label{well-definedness of M(x,n)}
Let $x$ be an indeterminate and $F(x)$ the free $\ld$-ring generated by $x$.
For any integer $m$ and any positive integer $n$, we have that
$\textsf{M}^m(x,n)\in F(x)$.
\end{lem}

\pf
Recall that
\begin{equation}\label{recursion of q necklace poly 0}
\begin{aligned}
n\textsf{M}^m(x,n)=-\sum_{d|n \atop d<n} d \left[\frac nd\right]_m \psi^{\frac nd}(\textsf{M}^m(x,d))+\left[n\right]_mx^n.
\end{aligned}
\end{equation}
First, assume that $1-m\ne 0.$
To show that it is divisible by $n$, choose an arbitrary prime divisor $p$ of $n$.
Then
\begin{align*}
n(1-m)\textsf{M}^m(x,n)
&\equiv -\sum_{d|\frac np} d(1-m^{\frac {n}{pd}})\psi^{\frac nd}(\textsf{M}^m(x,d))+(1-m^{\frac np})x^n\\
&\equiv -\psi^p (\sum_{d|\frac np} d(1-m^{\frac {n}{pd}})\psi^{\frac {n}{pd}}(\textsf{M}^m(x,d)))+(1-m^{\frac np})x^n \\
&\equiv -\psi^p ((1-m^{\frac np})x^{\frac np}))+(1-m^{\frac np})x^n \\
&\equiv (1-m^{\frac np})(x^n-\psi^p (x^{\frac np}))\\
&\equiv 0 \pmod{p^{\nu_p(n)}}.
\end{align*}
The last equivalence follows from the fact that
$$\textsf{M}^0(x,p^{\nu_p(n)})=\frac {1}{p^{\nu_p(n)}}(x^p-\psi^p(x^{\frac np}))\in F(x) \text{ (by Lemma \ref{to define Teichmuller map}}).$$
Consequently we have shown that $n(1-m)\textsf{M}^m(x,n)$ is divisible by $n$ and $1-m$.
To see that it is divisible by $n(1-m)$, let $\bar p$ be a common prime divisor of $n$ and $1-m$.
Let $n=\bar p^\nu n',\, d=\bar p^{\alpha}d'$ and $1-m=-\bar p^lk$, where $(n',\bar p)=(d',\bar p)=(k,\bar p)=1.$
Since
$$(1-m^{\frac{n}{d}})=1-(1+\bar p^lk)^{\bar p^{\nu-\alpha }\frac {n'}{d'}}\equiv 0 \pmod{\bar p^{\nu-\alpha+l}},$$
the product of the right hand side of Eq.\eqref{recursion of q necklace poly 0}) and $1-m$ is divisible by $p^{\nu+\alpha}$.
Hence, $n(1-m)\textsf{M}^m(x,n)$ is divisible by $p^{\nu+\alpha}$ and thus we can conclude that it is divisible by $n(1-m).$

If $m=1$, then
$$\sum_{d|n} \psi^{\frac nd}(\textsf{M}^1(x,d))=x^n.$$
By M\"{o}bius inversion formula, we get that
$$
\textsf{M}^1(x,n)=\sum_{d|n} \mu(d)\psi^d(x^{\frac nd}),\,\, (n\ge 1),
$$
which are obviously contained in $F(x)$
(see Lemma \ref{to define Teichmuller map}).
This completes the proof.
\qed
\vskip 2mm

For any $\ld$-ring $R$ and $r\in R$,
let $\pi: F(x) \to R$ be the $\ld$-ring homomorphism determined by
the condition $x \mapsto r$.
Define $\textsf{M}^m(r,n)$ by $\pi(\textsf{M}^m(x,n))$.
Consider the map
$\tau^m_R: W^m(R) \to \textsf B^m(R)$ defined by
$$(a_n)_{n\in \mathbb N}\mapsto \left(\sum_{d|n}\textsf{M}^m(a_d, \frac nd )\right)_{n\in \mathbb N}.$$
With the above preparation, we can prove the following theorem in the same way as in Theorem \ref{relation of W(g(q)) and B(g(q))}.

\begin{thm} \label{connection 1}
For each integer $m$, the assignment $R\mapsto \tau^m_R$ induces a natural isomorphism
$$\tau^m: W^m \circ F^{\textsf{$\ld$-Rings}}_{\textsf{Rings}} \to {\textsf B}^m|_{\textsf {$\ld$-Rings} }.$$
\end{thm}




\end{document}